\documentclass[12pt,reqno]{amsart}

\usepackage[T1]{fontenc}
\usepackage[utf8]{inputenc}
\usepackage{amsmath,amsfonts,amsthm,amssymb,amsxtra,bbm}
\usepackage{fourier,setspace,graphicx,color,pdflscape}
\usepackage[colorlinks=true]{hyperref}
\hypersetup{urlcolor=blue, linkcolor=blue, citecolor=red, anchorcolor=blue]}
\usepackage{todonotes}
\usepackage{fullwidth,mathrsfs}
\usepackage{amsthm}

\allowdisplaybreaks

\newtheorem{theorem}{Theorem}
\newtheorem{proposition}[theorem]{Proposition}
\newtheorem{lemma}[theorem]{Lemma}
\newtheorem{corollary}[theorem]{Corollary}
\newtheorem{remark}[theorem]{Remark}

\newtheorem{assumption}{Assumption}
\newcommand{\R}{{\mathbb R}}

\newcommand{\be}[1]{\begin{equation}\label{#1}}
\newcommand{\ee}{\end{equation}}

\newcommand{\ZDCT}{\mathrm{DC}}
\newcommand{\DC}{\mathrm{DC}}
\renewcommand{\leq}{\leqslant}
\renewcommand{\geq}{\geqslant}
\newcommand{\Id}{\mathrm{Id}}
\newcommand{\Lop}{\mathcal{L}}
\newcommand{\cF}{\mathcal{F}}
\newcommand{\rme}{\mathrm{e}}
\makeatletter
\@namedef{subjclassname@2020}{%
 \textup{2020} Mathematics Subject Classification}
\makeatother
\title[Decay rates kinetic Fokker--Planck]{How to construct explicit decay rates for kinetic Fokker--Planck equations?}

\author[Giovanni~Brigati]{Giovanni Brigati}
\address{\hspace*{-12pt}Giovanni ~Brigati: CEREMADE (CNRS UMR n$^\circ$7534), PSL university, Universit\'e Paris-Dauphine, Place de Lattre de Tassigny, 75775 Paris 16, France}
\email{brigati@ceremade.dauphine.fr}

\author[Gabriel Stoltz]{Gabriel Stoltz}
\address{\hspace*{-12pt}Gabriel Stoltz: CERMICS, Ecole des Ponts, Marne-la-Vallée, France \& MATHERIALS team-project, Inria Paris, France}
\email{gabriel.stoltz@enpc.fr}

\begin{document}
%%%%%%%%%%%%%%%%%%%%%%%%%%%%%%%%%%%%%%%%%%%%%%%%%%%%%%%%%%%%%%%%%%%%%%%%%%

\begin{abstract}
We study time averages for the norm of solutions to kinetic Fokker--Planck equations associated with general Hamiltonians. We provide fully explicit and constructive decay estimates for systems subject to a confining potential, allowing fat-tail, sub-exponential and (super-)exponential local equilibria, which also include the classic Maxwellian case. The key step in our estimates is a modified Poincar\'e inequality, obtained via a Lions--Poincar\'e inequality and an averaging lemma.
\end{abstract}

\maketitle
%%%%%%%%%%%%%%%%%%%%%%%%%%%%%%%%%%%%%%%%%%%%%%%%%%%%%%%%%%%%%%%%%%%%%%%%%%

\section{Introduction}\label{intro}

This work presents an extension of the hypocoercive approach studied in~\cite{brigati2022time} to the case of a non-zero potential energy function, and also considers the case of general kinetic energies beyond the standard quadratic one. Such generalizations can be relevant to enhance the performance of certain sampling methods in molecular dynamics and Bayesian statistics; see for instance~\cite{ST18,LFR19}. The more general framework also allows to make more transparent the structural assumptions behind the algebraic computations of~\cite{cao2023explicit,albritton2024variational}.

Kinetic models describe interacting multi-agent or multi-particle systems at an intermediate level between full microscopic and macroscopic scales \cite{balian2007microphysics}. More precisely, a kinetic partial differential equation (PDE) encodes the time evolution of the \textit{distribution function} of particles with respect to \textit{position} and some \textit{velocity} variable. Since the XIXth century, kinetic PDEs are formulated as a coupling between a transport and a collision/relaxation term. This is the case of the celebrated Boltzmann equation, for example. A general theory for Boltzmann's equation is still open. However, due to its physical validity, various simplified versions of it have been studied. A remarkable one is the Vlasov--Fokker--Planck equation, derived from Boltzmann's equation in a limiting regime corresponding to grazing collisions \cite{degond1992fokker}. Vlasov--Fokker--Planck equations date back at least to \cite{kolmogoroff1934zufallige}, and they are useful in statistical physics \cite{balian2007microphysics} and stellar dynamics~\cite{RevModPhys.15.1}. In general, Vlasov--Fokker--Planck equations should be understood as models for the time evolution of a distribution of particles subject to a potential and a random background force inducing relaxation towards equilibrium. This paper is concerned with a rather general class of Vlasov--Fokker--Planck equations introduced below.

\subsubsection*{Main convergence result}
We consider systems described by their configuration~$(x,v)$, where~$x \in \mathcal{X}$ is a spatial variable representing the positions of the particles composing the system, and~$v \in \mathcal{V}$ their velocities. Typical choices for the position space are $\mathcal{X} = \mathbb R^d$ for unconfined systems on the full space, or $\mathcal{X} = (L\mathbb{T})^d$ (with~$\mathbb{T} = \mathbb{R} \backslash \mathbb{Z}$ the unit torus) for systems enclosed in a cubic box of size~$L>0$ with periodic boundary conditions, a situation typical in molecular dynamics. The velocity space usually is~$\mathcal{V} = \mathbb{R}^d$, but here as well one could take~$\mathcal{V} = \mathbb{T}^d$, as considered for instance in \cite{ST18}. In the latter case, and more generally when considering non-quadratic kinetic energies, one should rather think of~$v$ as being a momentum, the variable conjugated to the positions, rather than a velocity. The Hamiltonian associated with the system is assumed to be separable, and is therefore the sum~$\phi(x) + \psi(v)$ of the potential energy function~$\phi:\mathcal{X}\to \mathbb{R}$ and the kinetic energy~$\psi:\mathcal{V}\to\mathbb{R}$. 
We introduce the following measures:
\[
\mu(dx) := \mathrm{e}^{-\phi(x)} \, dx, \qquad \gamma(dv) = \mathrm{e}^{-\psi(v)} \, dv.
\]
In order to have an invariant Boltzmann--Gibbs probability measure
\[
\Theta (dx \, dv) = \mu(dx) \, \gamma(dv),
\]
we shall assume that~$\mathrm{e}^{-\phi} \in \mathrm{L}^1(\mathcal{X})$ and~$\mathrm{e}^{-\psi} \in \mathrm{L}^1(\mathcal{V})$, and, without loss of generality, that $\mathrm{e}^{-\phi}$ and~$\mathrm{e}^{-\psi}$ integrate to~1, upon adding a constant to~$\phi,\psi$. The Vlasov--Fokker--Planck equation we study evolves a distribution~$f(t,x,v)$ according to
\begin{equation}
    \label{VFP}\tag{VFP}
    \left\{ \begin{aligned}
        & \partial_t f + \nabla_v \psi \cdot \nabla_x f - \nabla_x \phi \cdot \nabla_v f = \xi \, \nabla_v \cdot \left( \gamma \, \nabla_v \, \left( \frac{f}{\gamma} \right) \right), \\
        & f(t=0,\cdot,\cdot) = f_0 \in \mathrm{L}^2(\Theta^{-1}),
    \end{aligned} \right. 
\end{equation}
where $\xi > 0$ is a parameter regulating the intensity of the collision/relaxation term, and
\[
\mathrm{L}^2(\Theta) := \left\{ f \textrm{ measurable} \, \middle| \, \int_{\mathcal X
\times \mathcal V} |f|^2 \Theta < +\infty \right\}.
\]
We call 
$$Tf := \nabla_v \psi \cdot \nabla_x f  - \nabla_x \phi \cdot \nabla_v f, \qquad \mathscr Sf := \nabla_v \cdot \left( \gamma \, \nabla_v \, \left( \frac{f}{\gamma} \right) \right),$$
the transport and the collision operators, respectively. 

It is in fact more convenient to perform a change of unknown function and write~$f = h \Theta$, in which case we expect~$h$ to converge in the longtime limit to the constant function with value~$\int_{\mathcal{X} \times \mathcal{V}} f_0(x,v) \,dx \, dv$. The function~$h$ satisfies the following kinetic Ornstein--Uhlenbeck equation: 
\begin{equation}\label{VOU}\tag{VOU}
  \left\{ \begin{aligned}
    & \partial_t h + \nabla_v \psi \cdot \nabla_x h - \nabla_x \phi \cdot \nabla_v h = - \, \xi \,  \nabla_v^\star \, \nabla_v h, \\
    & h(t=0,\cdot,\cdot) = h_0 \in \mathrm{L}^2(\Theta),
  \end{aligned} \right.
\end{equation}
where~$A^\star$ denotes the adjoint of a closed operator on~$\mathrm{L}^2(\Theta)$, see \eqref{VOUo}. Equation \eqref{VOU} is linear and mass-preserving, hence we can assume the initial data to be of average zero with respect to $\Theta$ and study the long-time behaviour of the corresponding solution to \eqref{VOU}. 
We expect 
$$\lim_{t \to \infty} \|h(t,\cdot,\cdot)\|^2_{\mathrm{L}^2(\Theta)} = 0.$$
The goal of this work is to make this limit quantitative by obtaining effective decay rates. We sketch the typical result we obtain below (for actual statements, with general Hamiltonians and precise assumptions, see Section \ref{sec:main_results}).

\begin{theorem}[Informal statement]
    Let $\psi(v) = \frac{|v|^2}{2}$ and let $\phi$ be a smooth spatial potential, which grows linearly or faster as $|x| \to \infty$. Then, for all $\tau>0$, there exists an explicit constant $\lambda(\tau)>0$ such that all solutions to \eqref{VOU} satisfy
    \begin{equation*}
    \forall t > 0, \qquad \int_t^{t+\tau} \|h(s,\cdot,\cdot)\|^2_{\mathrm{L}^2(\Theta)} \, \frac{ds}{\tau} \, \leq \|h_0\|^2_{\mathrm{L}^2(\Theta)}   \, \mathrm{e}^{-\lambda(\tau) t}. 
    \end{equation*}
\end{theorem}

\subsubsection*{Relationship to other works}
In the recent literature, the problem of making the convergence rates of kinetic equations quantitative has been intensively investigated.
Equations such as \eqref{VOU} fall inside the abstract theory of \cite{hormander1967hypoelliptic}. From this founding paper, a recent theory originated, specialising on decay rates rather than regularity issues. This approach goes under the name of hypocoercivity. We briefly review the related literature, referring to the introduction of  \cite{bernard2020hypocoercivity} for a much more detailed account. Quantitative decay rates for some degenerate (but still hypoelliptic) diffusion equations appeared in~\cite{talay02,villani2006hypocoercive,villani2009hypocoercivity}. The strategy there relies on twisting the $\mathrm{H}^1$-norm in order to get an equivalent Lyapunov functional showing exponential decay in time along the flow of the equations under investigation. Then, exponential decay (at the same rate) can be recovered for the $\mathrm{H}^1$-norm, with a pre-factor $C>1$ arising from the twist between the reference norm and the Lyapunov functional.  

Later, a $\mathrm{L}^2-$setting for hypocoercivity was established \cite{MR2215889,dolbeault2015hypocoercivity}. The strategy is to perturb the standard $\mathrm{L}^2$-norm with a well-suited term in order to find an equivalent norm with a better dissipation estimate along kinetic PDEs. In addition, in this setting, one can study consistency with hydrodynamic limits and non-regularising kinetic PDEs with general initial data. The $\mathrm{L}^2$-framework is general enough to treat a large class of models~\cite{bouin2017hypocoercivity,bouin2019hypocoercivity,bouin2019fractional}, and it can provide sharp estimates for some specific cases \cite{arnold2014sharp,achleitner2016linear,achleitner2017optimal,achleitner2017multi,arnold2021}. 
Precise algebraic features of the framework and estimates with respect to the parameters defining the dynamics were studied in \cite{bernard2020hypocoercivity}. In particular, the scaling with respect to $\xi$ has been investigated in \cite{GROTHAUS20143515,dolbeault2013exponential,RS18,IOS19}. 
To conclude, in \cite{bouin20222} the $\mathrm{L}^2$-framework of hypocoercivity is analysed for a class of kinetic equations with factorised invariant measures, among which \eqref{VFP}, with specific choices of the potentials $\phi$ and $\psi$.

In order to complete this short review, let us finally mention a few other works related to hypocoercive techniques, namely \cite{eberle2019couplings, dietert2015convergence} for an approach based on couplings, \cite{baudoin2013bakryp} for a version of the Bakry--\'Emery methods for hypocoercive dynamics, and \cite{dietert2022quantitative} which recently studied a model with further degeneracy.   

\subsubsection*{Contributions of this work}
In the present paper, we focus on the approach of \cite{albritton2024variational}, which works directly in the standard $\mathrm{L}^2$-norm, without changing the reference norm as in, e.g., \cite{bouin20222}. This is possible thanks to an adapted space-time-velocity Poincar\'e inequality. We notice that time-integrated functionals of the solutions to kinetic equations are in use since \cite{guo2002landau, strain2004stability} at least. 
The framework of \cite{albritton2024variational} was qualitative although, constructive estimates can be derived in special cases \cite{brigati2022time,cao2023explicit}. One idea common to both references is to exploit a version of a space-time Poincar\'e--Lions inequality~\cite{amrouche2015lemma,carrapatoso2022weighted,dietert2022quantitative}, coupled with an explicit version of the computations in \cite[Proposition 6.2]{albritton2024variational}. 

In this paper we further generalise the techniques of \cite{brigati2022time,cao2023explicit}, by carefully keeping all estimates fully constructive. 
The main original contributions are:
\begin{enumerate}
    \item Adapting the Armstrong--Mourrat technique, we get fully explicit convergence estimates for \eqref{VOU} with general potential and kinetic energies $\phi,\psi$. In \cite{brigati2022time} the spatial variable was confined to a torus with $\phi=0,$ while \cite{cao2023explicit} only considers quadratic kinetic energies. Decay rates are shown to be exponential if $\psi$ grows fast enough at infinity. Otherwise, algebraic decay rates are obtained. 
    \item The potentials $\phi$ we consider are rather general. We remove the condition that the embedding $\mathrm{H}^1(\mu) \hookrightarrow \mathrm{L}^2(\mu)$ is compact, see \cite[Assumption 3]{cao2023explicit}. In addition, our decay rates are completely explicit, unlike those achieved via other approaches.
    \item Our estimates for the constants in the Lions inequality of Theorem \ref{thm2} should be compared with those of \cite{dietert2022quantitative}. Because performed in a simpler geometry, our result is fully constructive, as we actually build a solution for the elliptic equation involved to obtain the result. Moreover, as we avoid estimates of singular integrals \emph{\`a la} Calderon-Zygmund we obtain a better dimensional dependence in the constants. 
    \item The decay estimates are fully explicit in the fluctuation-collision parameter $\xi$. In particular, we discuss in Section~\ref{sec:discussinon_scaling_friction} that the convergence rate is of order~$\xi$ as~$\xi \to 0$, and of order~$\xi^{-1-1/\sigma}$ for~$\xi \to +\infty$, with~$\sigma$ related to the (heavy) tails of the probability distribution on the momenta.  
    \item For separable kinetic energies $\psi$, all constants are precisely quantified in terms of the dimension~$d$ and other parameters of the dynamics, generalising the analysis of~\cite{bernard2020hypocoercivity} and~\cite{cao2023explicit}.
\end{enumerate}
Some key technical points of our method are the following.
\begin{itemize}
\item We use time averages of the $\mathrm{L}^2$-norm of the solutions to \eqref{VOU} as entropies. Time averages allow for an effective entropy-entropy production estimate along the flow of \eqref{VOU}, where the key step in the proof is an adapted time-space-velocity Poincar\'e inequality (see Proposition~\ref{prop:Poincare_type}).
\item If $\psi$ grows at least linearly at infinity, the entropy production is strong enough to use a Gronwall lemma. Otherwise, the dissipation controls just a fractional power of the time average of the norm. Then, an algebraic decay is obtained via a Bihari--Lasalle argument. 
\item A fundamental ingredient is a space-time weighted Lions' inequality (see Theorem~\ref{thm2}), whose proof is fully constructive. Indeed, the problem is reduced to a regularity estimate for some elliptic PDE, which is solved either by explicit computations or with the Lax--Milgram theorem, depending on the nature of the source term.  
\item Our averaging Lemma \ref{avglem} extends and clarifies the algebraic structure behind analogous results in \cite{albritton2024variational,cao2023explicit}, working also for non-tensorised kinetic energies $\psi$. In particular, we provide a scheme to re-introduce the velocity dependence of a function to control the gradient of its space-time marginal. 
\end{itemize}

\subsubsection*{Extensions and perspectives} As mentioned earlier, our aim in this work is to present a general methodology to obtain decay rates for Fokker--Planck type equations which are fully explicit in the various parameters of the dynamics, in particular the friction and the dimension. Our presentation is as general as possible, so that extensions and adaptations to various other dynamics can be worked out. We have in mind here the linear Boltzmann equation, the linearized Landau equation~\cite{Rachid20}, Fokker--Planck equations with fractional Laplacians, generalized Langevin dynamics~\cite{OP11}, Adaptive Langevin dynamics~\cite{LSS19}, nonequilibrium dynamics~\cite{Dietert23}, etc. In fact, as our approach makes almost no use of the linearity of the Fokker--Planck equation at hand, generalizations to nonlinear equations can be envisioned.

\subsubsection*{Outline of the paper} This work is organized as follows. We start in Section~\ref{sec:main_results} by presenting the problem we consider and stating the main results we establish. We give in Section~\ref{sec:preliminary} some preliminary results which are useful for the proofs of the main results. Section~\ref{sec:sec3} is devoted to the proof of exponential decay estimates when the kinetic energy $\psi$ grows linearly or faster at infinity. Section~\ref{sec:sec4} proves algebraic decay rates for solutions to \eqref{VOU} for kinetic energies $\psi$ which grow slower than linearly at infinity. Weighted Lions estimates are proved in  Section \ref{sec:sec5}. Finally, a careful analysis of the scaling of the constants with respect of the dimension is given in Section \ref{sec:tenso} for tensorised kinetic energies. 

%------------------------------------------
\section{Main results}
\label{sec:main_results}

We give the main assumptions (Subsection~\ref{mr:ss1}) and results (Subsections~\ref{mr:ss15} to~\ref{mr:ss4}) of the paper. Finally, in Subsection~\ref{mr:ss5}, we discuss the dependence of our estimates on the parameters of the problem.

\begin{remark}
  The analysis we perform can be straightforwardly extended, at least from an algebraic viewpoint, to other types of degenerate symmetric dissipation operators, such as the one appearing in the linear Boltzmann equation, or involving fractional Laplacians, see \cite{dolbeault2015hypocoercivity,bouin2019fractional}.
\end{remark}

\subsection{Structural assumptions}\label{mr:ss1}
To present our main convergence result, we introduce a time~$\tau > 0$, and the associated uniform probability measure on the time interval~$[0,\tau]$:
\[
\mathrm{U}_\tau(dt) = \frac{1}{\tau} \mathbf{1}_{[0,\tau]}(t) \, dt.
\]
With some abuse of notation we sometimes denote probability measures and their densities with the same symbol. We also denote by $\mathrm{H}^{-1}(\mathrm{U}_\tau \otimes \mu)$ the dual of the space
\begin{equation}\label{dcdef}
\mathrm{H}^1_\ZDCT(\mathrm{U}_\tau \otimes \mu) = \left\{ h \in \mathrm{H}^{1}(\mathrm{U}_\tau \otimes \mu) \, \middle| \, h(0,\cdot) = h(\tau,\cdot) =0 \right\},
\end{equation}
where the subscript $\mathrm{DC}$ stands for \textit{Dirichlet conditions}. Note that the boundary conditions in time are well defined because the trace of functions~$h \in \mathrm{H}^{1}(\mathrm{U}_\tau \otimes \mu)$ makes sense on the boundary~$\{0,\tau\} \times \mathcal{X}$ of the domain. We state the main assumptions we need for the potential energy~$\phi$ and the kinetic energy~$\psi$. 

\begin{assumption}\label{ass1}
  The function $\phi$ is smooth and $\mathrm{e}^{-\phi} \in \mathrm{L}^1(dx)$, with $\int_{\mathcal{X}} \mathrm{e}^{-\phi(x)} \, dx =1$. Moreover, for any time $\tau>0$, there exists a constant~$C^\mathrm{Lions}_\tau>0$ such that the following Lions inequality holds true:
 \begin{equation}
    \label{lionstx}
    \forall g \in \mathrm{L}^2(\mathrm{U}_\tau \otimes \mu), \quad \left\|g- \iint_{[0,\tau] \times \mathcal{X}} g(t,x) \, \mathrm{U}_\tau(dt) \, \mu(dx)\right\|^2_{\mathrm{L}^2(\mathrm{U}_\tau \otimes \mu)} \leq C_\tau^\mathrm{Lions} \|\nabla_{t,x} g  \|^2_{\mathrm{H}^{-1}(\mathrm{U}_\tau \otimes \mu)}.
\end{equation}
\end{assumption}

Sufficient conditions for \eqref{lionstx} to  hold with explicit estimates of $C^{\mathrm{Lions}}_\tau$ are discussed in Section \ref{mr:ss4}.

\begin{assumption}\label{ass2}
The function~$\psi$ is smooth and $\mathrm{e}^{-\psi} \in \mathrm{L}^1(dv)$, with $\int_{\mathcal{V}} \mathrm{e}^{-\psi(v)} \, dv =1$. Moreover,
\begin{equation}
    \label{eq:integrability_conditions_on_psi}
    \int_{\mathcal{V}} |\nabla_v \psi(v)|^4 \, \gamma(dv) +
    \int_{\mathcal{V}} \left|\nabla_v^2 \psi(v)\right|^2 \, \gamma(dv) < +\infty,
    \qquad
    \limsup_{R \to \infty} \, \int_{\mathcal{V} \setminus B(0,R)} |\nabla \psi|^{-2} \, \gamma(dv) = 0,
  \end{equation}
  and the symmetric positive matrix
  \begin{equation}
    \label{eq:matrix_M}
    \mathscr{M} = \int_{\mathcal{V}} \nabla_v \psi \otimes \nabla_v \psi \, d\gamma
  \end{equation}
is definite.
\end{assumption}

The condition~\eqref{eq:matrix_M} is satisfied when~$\mathcal{V} = \mathbb{R}^d$, see~\cite[Appendix~A]{ST18} for a proof. More generally, this condition holds as soon as there is no direction~$\rho \in \mathbb{R}^d$ such that~$\rho^\top \nabla\psi(v) = 0$ for all~$v \in \mathcal{V}$ (see the argument in the proof of~\cite[Theorem~3.3]{bernard2020hypocoercivity}). The integrability conditions~\eqref{eq:integrability_conditions_on_psi} are rather easy to satisfy. The estimates hold for instance if~$\mathcal{V}$ is bounded, or~$\psi$ increases polynomially or logarithmically at infinity provided it grows fast enough; see~\eqref{exmp1} and~\eqref{exmp2} below for more precise statements.

\begin{assumption}\label{ass3}
The Lie algebra generated by $\nabla_v$ and $\nabla_v \psi \cdot \nabla_x  - \nabla_x \phi \cdot \nabla_v$ is of dimension~$2d$ at all points $(x,v) \in \mathcal{X} \times \mathcal{V}$.
\end{assumption}

Sufficient conditions ensuring that the Lie algebra is full at every point are given in Appendix~\ref{app:Hormander} for completeness.

\subsection{An averaging lemma}\label{mr:ss15}
Averaging lemmas are a typical tool in kinetic theory to recover increased regularity on the velocity average of the distribution function. This is the idea behind our first result, which is fundamental for the main theorems in Sections~\ref{mr:ss2} and~\ref{mr:ss3}. 

We need a further assumption on $\phi$. 

\begin{assumption}\label{ass1bis}
There exists a constant $L_\phi \in \mathbb R_+$ such that
\begin{equation}\label{commphi}
    \forall z \in \mathrm{H}^1(\mu), \qquad \|z \, \nabla_x \phi \|_{\mathrm{L}^2(\mu)} \leq L_\phi \|z\|_{\mathrm{H}^1(\mu)}.
\end{equation}
\end{assumption}
The inequality \eqref{commphi} is satisfied (with explicit control on the constant) in many cases of interest, see ~\cite[Lemma~A.24]{villani2009hypocoercivity},  as well as~\cite[Lemma~2.2]{cao2023explicit} and~\cite[Lemma~3.7]{bernard2020hypocoercivity} for a precise quantification of the dependence of the constant $L_\phi$ on the dimension, under certain structural assumptions on the potential $\phi$. This is recalled in Lemma~\ref{lemmbfls} below.

A functional space which will be useful in the proofs is the kinetic space, defined as follows for a positive time $\tau>0$:
\[
\mathrm{H}^1_{\mathrm{kin}} = \left\{ h \in \mathrm{L}^2\left(\mathrm{U}_\tau \otimes \mu, \mathrm{H}^1(\gamma)\right) \ \middle| \ (\partial_t + T)h \in \mathrm{L}^2\left(\mathrm{U}_\tau \otimes \mu, \mathrm{H}^{-1}(\gamma)\right) \right\}.
\]
We also introduce the following projector on~$\mathrm{L}^2(\mathrm{U}_\tau \otimes \Theta)$, corresponding to velocity averaging:
\[
\Pi h = \int_{\mathcal{V}} h(\cdot,\cdot,v)  \,\gamma(dv).
\]
We finally denote by $\mathrm{L}^2_0(\mathrm{U}_\tau \otimes \Theta)$ the space of functions in $\mathrm{L}^2(\mathrm{U}_\tau \otimes \Theta)$ with average~0 with respect to~$\mathrm{U}_\tau \otimes \Theta$. We are then in position to state the following useful averaging lemma for~\eqref{VOU}.

\begin{lemma}
  \label{avglem}
  Under Assumptions~\ref{ass1}, \ref{ass2} and~\ref{ass1bis}, there exists an explicit constant $K_\mathrm{avg} \in \mathbb{R}_+$ such that, for any~$h \in \mathrm{H}^1_{\mathrm{kin}}$,
  \[
  \|\nabla_{t,x} \Pi\, h\|^2_{\mathrm{H}^{-1}(\mathrm{U}_\tau \otimes \mu)} \leq K_\mathrm{avg} \left( \|(\Id-\Pi)h\|^2_{\mathrm{L}^2(\mathrm{U}_\tau \otimes \Theta)} + \|(\partial_t+T)h \|^2_{\mathrm{L}^2(\mathrm{U}_\tau \otimes \mu, \mathrm{H}^{-1}(\gamma))} \right).
  \]
\end{lemma}

The proof of this result can be read in Section~\ref{sec:avglem}. In particular, we provide in~\eqref{kavg} a general explicit upper bound on $K_{\mathrm{avg}}$. A precise dependence in terms of the dimension can be obtained for separable kinetic energies, as discussed in Section~\ref{sec:scaling_dimension}; see Proposition~\ref{kavgtens} in Section~\ref{sec:tenso}.

\subsection{Exponential convergence rates} \label{mr:ss2}

We consider the longtime behaviour of solutions to \eqref{VOU}. In particular, our main goal is to prove that 
$$h(t,\cdot,\cdot) \xrightarrow[t \to \infty]{} \iint_{\mathcal{X} \times \mathcal{V}} h_0(x,v) \, \Theta(dx\,dv),$$
for all solutions to \eqref{VOU}, with quantitative and explicit convergence rates. 
As explained in Section \ref{intro}, this reduces to studying the decay to $0$ for solutions to~\eqref{VOU} with initial condition
\[
h_0 \in \mathrm{L}^2_0(\Theta) = \left\{ g \in \mathrm{L}^2(\Theta) \ \middle| \ \int_{\mathcal{X} \times \mathcal{V}} g \, d\Theta = 0\right\}.
\]
Other sub-spaces such as $\mathrm{L}^2_0(\gamma)$ are defined in a similar way.
In order to achieve exponential decay rates, we need the following. 
\begin{assumption}\label{ass4}
A Poincar\'e inequality holds for~$\gamma$: there exists $c_\psi>0$ such that 
\begin{equation}
  \label{eq:coercivity_Delta_psi_OU}
  \forall g \in \mathrm{H}^1(\gamma) \cap \mathrm{L}^2_0(\gamma), \qquad \|g\|_{\mathrm{L}^2(\gamma)} \leq c_\psi^{-1/2} \|\nabla_v g \|_{\mathrm{L}^2(\gamma)}.
\end{equation}
\end{assumption}
A sufficient condition on~$\psi$ for the latter inequality to hold is that there exists a constant $a \in (0,1)$ such that
$$ \liminf_{|v| \to \infty} \left( a|\nabla_v \psi(v)|^2 - \Delta_v \psi(v) \right) > 0,$$
see~\cite{persson1960bounds,10.1214/ECP.v13-1352}. The inequality~\eqref{eq:coercivity_Delta_psi_OU} is always true when~$\mathcal{V} = \mathbb{T}^d$ and~$\psi$ is smooth. 

As in \cite{brigati2022time}, following up on ideas used in \cite{guo2002landau,albritton2024variational}, we introduce, for any $h \in \mathrm{L}^2([0,\infty); \mathrm{L}^2(\Theta)),$ the following functional:
\begin{equation}\label{timeavg}
\mathscr{H}_\tau(t) = \int_t^{t+\tau} \|h(s,\cdot,\cdot)\|^2_{\mathrm{L}^2(\Theta)} ds.
\end{equation}
Then, our main result goes as follows (see Section~\ref{sec:proof_thm1} for the proof). 

\begin{theorem}
  \label{thm1new}
  Under Assumptions \ref{ass1}, \ref{ass2}, \ref{ass3}, \ref{ass1bis} and~\ref{ass4}, there exists a constant $\lambda>0$ such that all functions $h$ solving \eqref{VOU} satisfy  
  \[
  \forall \, t \geq 0, \qquad \mathscr{H}_\tau(t) \leq  \mathrm{e}^{-2\lambda \, t} \mathscr{H}_\tau(0).
  \]
  Moreover, 
  \begin{equation}
    \label{eq:bar_lambda}
    \lambda \geq  \bar{\lambda} := \frac{1}{1 + C^{\mathrm{Lions}}_\tau K_\mathrm{avg}}\left( \frac{1}{\xi c_\psi} + \xi \right)^{-1}.
  \end{equation}
\end{theorem}

The dependence of the right hand side of~\eqref{eq:bar_lambda} with respect to various key parameters is made precise in Section~\ref{mr:ss5}. This includes the scaling with respect to the dimension~$d$ for separable kinetic energies, and the scaling with respect to the friction~$\xi$ for general kinetic energies. 
The latter scaling is fulling determined by the last
factor on the right hand side of~(2.8) since the
constants~$K_\mathrm{avg}$ and~$C^\mathrm{Lions}_\tau$ appear in
functional inequalities involving functions of~$t,x$ only.

Theorem~\ref{thm1new} implies a standard hypocoercivity estimate (see \cite{villani2009hypocoercivity}) for solutions to \eqref{VOU}. Note that there is no prefactor~$C$ in the exponential decay estimate of Theorem~\ref{thm1new}. Pointwise in time decay estimates can be deduced from Theorem~\ref{thm1new} upon adding such a prefactor (see Section~\ref{subseccor1} for the proof).

\begin{corollary}\label{cor1new}
  Under the same assumptions as Theorem \ref{thm1new}, all solutions to \eqref{VOU} satisfy
  \[
  \forall t \geq 0, \qquad  \|h(t,\cdot,\cdot)\|^2_{\mathrm{L}^2(\Theta)} \leq \mathrm{e}^{\lambda \tau} \, \mathrm{e}^{-\lambda t} \, \|h_0\|^2_{\mathrm{L}^2(\Theta)}.
  \]
\end{corollary}

\subsection{Algebraic convergence rates}\label{mr:ss3}

Condition \eqref{eq:coercivity_Delta_psi_OU} is equivalent to the fact that the symmetric part~$\mathcal{S}$ of the Fokker--Planck operator is coercive, in the following sense:
\begin{equation}
\label{eq:coercivity_Delta_psi}
  \mathcal{S} \geq c_\psi \mathbb{1} \quad \text {on} \quad  \mathrm{L}^2_0(\gamma).
\end{equation}
Typical cases where the last condition (i.e. Assumption \ref{ass4}) does not hold correspond to a kinetic energy~$\psi$ growing less than linearly at infinity, such as
\begin{equation}
  \label{exmp1}
    \psi(v) = \left(1+|v|^2\right)^{\frac{\alpha}{2}}, \qquad \alpha \in (0,1), 
\end{equation}    
or    
\begin{equation}
  \label{exmp2}
  \psi(v) = \frac{\beta}{2} \, \log \left(1+|v|^2\right) , \qquad \beta >d+4, 
\end{equation}  
so that $\gamma$ is heavy tailed at infinity. The fact that~$\beta$ needs to be strictly larger than~$d+2$ in~\eqref{exmp2} comes from the requirement that~$\mathrm{e}^{-\psi}$ be integrable; and in fact we need~$\beta > d+4$ in order for the second condition in~\eqref{eq:integrability_conditions_on_psi} to hold. Although the two examples~\eqref{exmp1} and~\eqref{exmp2} do not lead to a probability measure~$\gamma$ satisfying a Poincar\'e inequality, the associated probability measures satisfy the following weighted Poincar\'e inequality (see \cite{blanchet2007hardy,bouin2019hypocoercivity}). 

\begin{assumption}\label{ass5}
There exist a constant $P_{\psi} \in \mathbb R_+$ and a function $\mathscr G \, \in \mathrm{L}^1(\gamma)$ such that~$\mathscr{G} \geq 1$ and 
\begin{equation}\label{eq:weightedpoi}
    \forall g \in \mathrm{H}^1(\gamma) \cap \mathrm{L}^2_0(\gamma), \qquad \int_{\mathbb R^d} \mathscr G^{-1}(v) \, g^2(v) \, \gamma(dv) \leq P_\psi \|\nabla_v g \|^2_{\mathrm{L}^2(\gamma)}.
\end{equation}
\end{assumption}

Note that, as emphasized in~\cite[Appendix~A]{bouin2019hypocoercivity} for instance, the function~$g$ in~\eqref{eq:weightedpoi} is not centered with respect to the probability measure proportional to~$\mathscr G^{-1} \, d\gamma$, but with respect to the probability measure~$d\gamma$ which appears on the right hand side.
For the classes of examples \eqref{exmp1} and \eqref{exmp2}, we can respectively consider $\mathscr G(v) = (1+|v|^2)^{1-\alpha}$ and $\mathscr G(v) = 1+|v|^2$, see~\cite{blanchet2007hardy,bouin2019hypocoercivity}.

In the cases covered by Assumption \ref{ass5}, we recover constructive algebraic decay rates for~$\mathscr H_\tau$ as follows.

\begin{theorem}\label{thm2new}
Suppose that Assumptions~\ref{ass1}, \ref{ass2}, \ref{ass3}, \ref{ass1bis} and~\ref{ass5} hold, and consider $\sigma > 0$ such that $\mathscr G^\sigma \in \mathrm{L}^1(\gamma)$. For all solutions~$h$ to~\eqref{VOU} with~$h_0 \in \mathrm{L}^\infty(\mathcal{X} \times \mathcal{V})$, there exist two explicit and positive constants $c_1(h_0),c_2(h_0)$ depending only on~$\|h_0\|_{\mathrm{L}^\infty(\mathcal{X} \times \mathcal{V})}, \phi, \psi$ and~$\sigma$ such that
\[
\forall \xi > 0, \quad \forall t \geq 0, \qquad   \mathscr{H}_\tau(t) \leq \frac{\mathscr{H}_\tau(0)}{\left( 1 + \left( \xi^{-\frac{\sigma}{\sigma+1}} c_1 + \xi c_2 \right)^{-\frac{\sigma+1}{\sigma}} t \right)^{\sigma}}.
\]
\end{theorem}

The proof of this result can be read in Section~\ref{sec:sec4}, where, in particular, explicit expressions for $c_1,c_2$ are given (see~\eqref{eq:expressions_c1_c2}). As in Section \ref{mr:ss3}, a decay estimate for $\mathscr{H}_\tau$ implies a pointwise decay estimate in time. We state a result analogous to the one of Corollary \ref{cor1new} in the present situation (see also Section~\ref{sec:sec4} for the proof).

\begin{corollary}
  \label{cor2new}
  Under the same assumptions as Theorem \ref{thm2new}, all solutions to \eqref{VOU} satisfy
  \[
  \forall t \geq \tau, \qquad \|h(t,\cdot,\cdot)\|^2_{\mathrm{L}^2(\Theta)} \leq \frac{\|h_0\|^2_{\mathrm{L}^2(\Theta)}}{\left( 1 + \left( \xi^{-\frac{\sigma}{\sigma+1}} c_1 + \xi c_2 \right)^{-\frac{\sigma+1}{\sigma}} (t-\tau) \right)^{\sigma}}.
  \]
\end{corollary}

\begin{remark}
  Note that Theorems~\ref{thm1new} and~\ref{thm2new} are qualitatively consistent. Indeed, the first one can be recovered from the second one in the regime corresponding to a sequence of sub-exponential local equilibria~$\gamma$ becoming exponential in the limit. In this case, \eqref{eq:weightedpoi} holds with $\mathscr G \to 1$ and~$\sigma \to +\infty$. The denominator of the term on the right hand side of the estimate in Corollary~\eqref{cor2new} then formally scales as $[(t-\tau)/(c_1 \xi^{-1} + c_2\xi)]^\sigma$, which suggests an exponential decay governed by a quantity of the form~$\min(\xi,\xi^{-1})t$.
\end{remark}

\begin{remark}
  An alternative approach to obtaining a convergence result similar to Theorem~\ref{thm2new} is that of weak Poincar\'e inequalities, as considered in~\cite{grothaus2019weak}; see~\cite{BSWW24} for an adaptation of the approach considered in this work to that case.
\end{remark}

Theorem~\ref{thm2new} crucially relies on the $\mathrm{L}^\infty$ bound on initial data. When the local equilibrium $\gamma$ decays slower than exponentially at infinity, the phenomenon of \textit{loss of velocity moments} occurs. This means that the $\mathrm{L}^2-$norm of the solution cannot be controlled without assuming the control of extra velocity moments for the initial data. In addition, we have to ensure that these moments are propagated by \eqref{VOU}. Such an issue is present even in the spatially-homogeneous case. This has already been studied in~\cite{bouin2019hypocoercivity,brigati2022time,carrapatoso2017landau}. In order to follow the strategy of proof of~\cite{bouin2019hypocoercivity,brigati2022time}, the term one needs to bound along the flow of~\eqref{VOU} is
\begin{equation}\label{momentumg}
  \iint_{\mathcal{X} \times \mathcal{V}} \mathscr G(v)^\sigma \, h^2(t,x,v) \, \Theta(dx \, dv).    
\end{equation}
However, obtaining a direct uniform control, which is also quantitative, requires substantial work. This is why, for simplicity, we restrict ourselves to bounded initial data and to apply the strong maximum principle for hypoelliptic operators. A more general approach, based on Lyapunov functionals and aimed at removing the~$\mathrm{L}^\infty$ assumption on the initial condition, is currently investigated in~\cite{brigati2025moments}.
%----------------------------
\subsection{Weighted Lions' inequalities}\label{mr:ss4}

The inequality~\eqref{lionstx} is a generalised version of Lions' inequality~\cite{amrouche2015lemma}. Such an inequality with a weighted probability measure has first been studied in~\cite{carrapatoso2022weighted}. We give here simple sufficient conditions on the potential~$\phi$ for~\eqref{lionstx} to hold. We also obtain an explicit estimate of $C^{\mathrm{Lions}}_\tau$, under some growth conditions on the spatial potential~$\phi$.

\begin{theorem}
  \label{thm2}
  Assume that~$\phi$ is a smooth potential energy function such that~$\mathrm{e}^{-\phi} \in \mathrm{L}^1(\mathcal{X})$ with
  \[
  \int_\mathcal{X} \mathrm{e}^{-\phi(x)} \,dx = 1,
  \]
  that its associated probability measure $\mu$ satisfies a Poincar\'e inequality, namely
  \begin{equation}
    \label{poiphi}
    \forall z \in \mathrm{H}^1(\mu) \cap \mathrm{L}_0^2(\mu), \qquad \|z\|_{\mathrm{L}^2(\mu)} \leq c_\phi^{-1/2}\|\nabla_x z\|_{\mathrm{L}^2(\mu)},
  \end{equation}
  for some $c_\phi > 0$, and that
  \begin{equation}\label{hess}
    |\nabla^2 \phi|^2 \leq (c'_\phi)^2 \, \left( \, d + |\nabla_x \phi|^2 \right), \qquad \Delta \phi \leq c_\phi'' \left( d +  |\nabla_x \phi|^2 \right),
  \end{equation}
  for some constants $c'_\phi, c''_\phi \in \R_+$.
  Then, Assumption \ref{ass1} is satisfied, with
  \[
  C_\tau^{\mathrm{Lions}} = C_\phi \left(\frac{1}{\tau^2} + \left(1+\tau^2\right) \sqrt{d} \right),
  \]
for a constructive constant $C_\phi$, depending only on $c_\phi,c'_\phi,c''_\phi$, but not directly on $d$.
\end{theorem}

A precise upper bound for~$C_\tau^{\mathrm{Lions}}$ is provided in Proposition~\ref{prop:estimates_divergence_equation}, see the formulae~\eqref{eq:final_constant_Lions_ineq} and~\eqref{eq:C_tau_LM}-\eqref{eq:C_tau_N}.

The fact that the Lions inequality holds is a direct consequence of the results obtained in the more general setting of~\cite{dietert2022quantitative}. Our main contribution here is to precisely track the dependence on the various parameters through an explicit construction of the solution, which is possible here given the simple geometry of the space-time domain we consider. 

The scaling of the constants in~\eqref{hess} is natural, as it is the correct one for separable potentials
\begin{equation}
  \label{eq:separable_pot}
  \phi(x) = \sum_{i=1}^d \bar{\phi}(x_i),
\end{equation}
with $\bar{\phi} : \R \to \R$ such that $\left|\bar{\phi}''\right| \leq c''_\phi \left( 1 + |\bar{\phi}'|\right)$. Sufficient conditions for~\eqref{poiphi} to hold true are given in~\cite{persson1960bounds,10.1214/ECP.v13-1352}, see the discussion after Assumption~\ref{ass4}.

The explicit dependence on $d$ in the estimate of $C^{\mathrm{Lions}}_\tau$ (which is exactly~$\sqrt{d}$ if the Poincar\'e constant $c_\phi$ does not depend on $d$) allows to clearly control the dimension dependence in the convergence rate of Theorem~\ref{thm1new}. Dimension independent Poincar\'e constants can be obtained when~$\mu$ is the tensorisation of $1$-dimensional measures, or for systems with finite number of interactions, see the discussion in \cite[Section 3.1.2]{bernard2020hypocoercivity}. A careful analysis of our estimates for~$C_\tau^\mathrm{Lions}$ shows that this quantity scales as~$1/c_\phi$ when~$c_\phi \to 0$ (see Remark~\ref{rmk:scaling_Lions_constant_Poincare}). 
  
\subsection{Conclusions}\label{mr:ss5}
We conclude the presentation of our results by briefly discussing how our estimates explicitly depend on some key parameters. 

\subsubsection{Scaling with respect to the dimension}
\label{sec:scaling_dimension}
In order to more carefully identify the scaling of the exponential convergence rate in Theorem~\ref{thm1new} with respect to the dimension~$d$, and provided the estimates~\eqref{poiphi}-\eqref{hess} hold with a known dimension dependence of the Poincar\'e constant~$c_\phi$, it suffices to study how~$K_\mathrm{avg}$ depends on the dimension. We restrict ourselves to this end to a simpler setting where~$\mathcal{V} = \mathbb R^d$. Our analysis could however be easily extended to more general spaces~$\mathcal{V}$, with suitable adaptations in the notation. More importantly, we assume that the kinetic energy~$\psi : \mathcal{V} \to \mathbb R$ is separable, \emph{i.e.} 
\begin{equation}
  \label{tensorised}
  \psi(v) = \sum_{i=1}^d q(v_i),
\end{equation}
where $q: \mathbb R \to \mathbb R$ is a real valued function. In this context, Assumption~\ref{ass2} is implied by the conditions
\[
\mathrm{e}^{-q} \in \mathrm{L}^1(\mathbb R), \qquad \int \mathrm{e}^{-q(v_1)} \, dv_1 =1,
\]
and
\begin{equation}
  \label{eq:integrability_conditions_on_psi2}
  \begin{aligned}
  & \int_{\mathbb R} \left( \left|q'(v_1)\right|^4 + \left|q''(v_1)\right|^2 \right) \, \mathrm{e}^{-q(v_1)} \, d v_1 < +\infty, \\
  & \limsup_{R \to \infty} \, \int_{(-\infty,-R) \, \cup \, (R,\infty)} \left|q'(v_1)\right|^{-2} \, \mathrm{e}^{-q(v_1)} \, dv_1= 0.
  \end{aligned}
\end{equation}
Moreover, Assumption~\ref{ass4} is equivalent to the fact that the following one-dimensional Poincar\'e inequality holds:
\begin{equation}
  \label{eq:coercivity_Delta_psi_OU2}
  \forall g \in \mathrm{H}^1\left(\mathrm{e}^{-q}\right) \cap \mathrm{L}^2_0\left(\mathrm{e}^{-q}\right), \qquad \|g\|_{\mathrm{L}^2(\mathrm{e}^{-q})} \leq c_\psi^{-1/2} \| g' \|_{\mathrm{L}^2(\mathrm{e}^{-q})}.
\end{equation}
Indeed, \eqref{eq:coercivity_Delta_psi_OU2} implies by tensorisation~\eqref{eq:coercivity_Delta_psi_OU} with the same constant $c_\psi$; while~\eqref{eq:coercivity_Delta_psi_OU} applied to a function of the first component~$v_1$ only leads to~\eqref{eq:coercivity_Delta_psi_OU2}. In particular, $c_\psi$ does not depend on the dimension~$d$.

\begin{proposition}
  Suppose that the kinetic energy is of the separable form~\eqref{tensorised} with the conditions~\eqref{eq:integrability_conditions_on_psi2} and~\eqref{eq:coercivity_Delta_psi_OU2} satisfied, and additionally that
  \begin{equation}
    \label{eq:condition_vanishing}
    \int_\mathbb{R} q^{(4)} \, \rme^{-q} = 0.
  \end{equation}
  Suppose moreover that the potential energy is smooth with~$\rme^{-\phi} \in \mathrm{L}^1(\mathcal{X})$, and such that~\eqref{poiphi} and~\eqref{hess} hold true. Suppose finally that Assumption~\ref{ass3} holds. Then, the lower bound~$\bar{\lambda}$ to the exponential decay rate in Theorem~\ref{thm1new} scales as 
  \[
  \bar{\lambda} = \frac{c_\phi}{1 + K_\mathrm{avg}^q R^\mathrm{Lions}}\left(\frac{1}{c_\psi \xi} + \xi\right),
  \]
  where~$K_\mathrm{avg}^q$, given in~\eqref{eq:K_avg_q}, is independent of the dimension~$d$, and
  \[
  R^\mathrm{Lions} = C_1\left(c_\phi',c_\phi''\right) \left(c_\phi+\sqrt{d}\right) + C_2,
  \]
  for some constants~$C_1,C_2$ which are also independent of the dimension~$d$ when~$c_\phi',c_\phi''$ are independent of~$d$; see~\eqref{eq:R_Lions}.
\end{proposition}

The proof of this result is a direct consequence of~\eqref{eq:bar_lambda}, Theorem~\ref{thm2} and precise estimates on~$K_{\mathrm{avg}}$ as given by Proposition~\ref{kavgtens} in Section~\ref{sec:tenso}. Note that we need to consider the additional assumption~\eqref{eq:condition_vanishing} in order to ensure that the constant~$K_{\mathrm{avg}}$ can be chosen independently of the dimension. When~\eqref{eq:condition_vanishing} is not satisfied, the constant~$K_{\mathrm{avg}}$ scales as~$\sqrt{d}$ (as made precise in Remark~\ref{rmk:dim_dep_Kavg}) and the dimension dependence of the convergence rate we obtain via Theorem~\ref{thm2} degrades. Let us emphasize here that we do not claim that the dimension dependence we get from our technique of proof is optimal, as it could potentially be beneficial to simultaneously perform the manipulations used to establish the averaging lemma and the Lions inequality, instead of doing them sequentially. This could indeed alleviate the need to estimate some terms, as in the proof of~\cite[Theorem~2]{cao2023explicit}.

When~\eqref{eq:condition_vanishing} is satisfied (which is the case for the most frequent choice in practice, namely a Gaussian distribution of velocities), the result we obtain is consistent with the one in~\cite{cao2023explicit}, with however fully explicit rates here. 

\subsubsection{Scaling with respect to the time average parameter} The dependence on $\tau$ in the constants appears only through~$C_\tau^{\mathrm{Lions}}$, which scales as~$\max(\tau^2,\tau^{-2})$. It is not possible to pass to the limits~$\tau \to 0$ or~$\tau \to +\infty$ in our estimates. For instance, $\bar{\lambda}$ in Theorem~\ref{thm1new} scales as~$\frac{1}{1+C_\tau^{\mathrm{Lions}}}$, and therefore goes to~0 as~$\tau \to 0$ or~$\tau \to +\infty$. In particular, the degeneracy arising from the limit~$\tau \to 0$ forbids a direct passage to a pointwise exponential decay with $C=1$ from Theorem~\ref{thm1new}, which would anyway lead to an absurd statement, as the evolution operator fails to be coercive.

Let us emphasize that we do not claim that the scaling of our bounds are sharp with respect to the time~$\tau$. However, we improve on \cite[Equation~(27)]{cao2023explicit} in the factor $\tau^{-2}$, achieving the conjectured scaling of~\cite[Remark 2.7]{cao2023explicit} in Proposition~\ref{prop:estimates_divergence_equation} below.  %Better estimates are obtained in~\cite{cao2023explicit}.{impli
%  We do not claim that the dependence in~$\tau$ in our estimates for~$C_\tau^{\mathrm{Lions}}$ are sharp. Our estimate cannot be compared to the one provided by~\cite[Equation~(27)]{cao2023explicit} as the time derivative of some of the terms is missing.

\subsubsection{Scaling with respect to the friction}
\label{sec:discussinon_scaling_friction}
In order to highlight the role of~$\xi$, we briefly discuss a microscopic derivation of \eqref{VFP}, seen as the PDE governing the evolution of the law of Langevin dynamics. The Langevin dynamics associated with the Hamiltonian~$\phi(x)+\psi(v)$ read
\begin{equation}
  \label{eq:Langevin}
\left\{ \begin{aligned}
  dX_t & = \nabla_v \psi(V_t) \, dt, \\
  dV_t & = -\nabla_x \phi(X_t) \, dt - \xi \nabla_v \psi(V_t) \, dt + \sqrt{2\xi} \, dW_t,
\end{aligned} \right.
\end{equation}
where~$\xi > 0$ is the friction coefficient and~$(W_t)_{t \geq 0}$ is a standard $d$-dimensional Brownian motion. Let us motivate that we expect the relaxation time of~\eqref{eq:Langevin} to the stationary state to be of order~$\max(\xi,\xi^{-1})$. In the limit~$\xi \to 0$, Langevin dynamics reduce to Hamiltonian dynamics. For small frictions, the average energy evolves on timescales of order~$\xi^{-1}$ since, by Itô calculus,
\[
\frac{d}{dt} \mathbb{E}\left[H(X_t,V_t)\right] = \mathbb{E}\left[(T+\xi\mathcal{S})H)(X_t,V_t)\right] = \xi \mathbb{E}\left[(\mathcal{S} H)(X_t,V_t)\right],
\]
where~$T + \xi \mathcal{S}$ is the generator of~\eqref{eq:Langevin}. The other limit of interest is the overdamped regime~$\xi \to +\infty$, where one recovers the overdamped Langevin dynamics upon rescaling time by a factor~$\xi$:
\[
d\overline{X}_t = -\nabla_x \phi(\overline{X}_t) \, dt + \sqrt{2} \, dB_t,
\]
where~$(B_t)_{t \geq 0}$ is also a standard $d$-dimensional Brownian motion. The formal proof of this limit relies on the observation that solutions of~\eqref{eq:Langevin} satisfy (upon replacing the term~$- \xi \nabla_v \psi (V_t) \, dt$ in the second line of~\eqref{eq:Langevin} by~$-\xi \,dX_t$ and integrating in time from~0 to~$\xi t$)
\[
\begin{aligned}
X_{\xi t} - X_0 &= \frac{V_0-V_{\xi t}}{\xi} - \frac{1}{\xi} \int_0^{\xi t} \nabla_x \phi(X_s) \, ds + \sqrt{\frac{2}{\xi}} \, dW_ {\xi s} \\
&= \frac{V_0 - V_{\xi \, t}}{\xi} - \int_0^t \nabla_x \phi (X_{\xi \,s}) ds + \sqrt{\frac{2}{\xi}} \, dW_ {\xi s},
\end{aligned}
\]
where~$(\xi^{-1/2} W_{\xi t})_{t \geq 0}$ still is a standard $d$-dimensional Brownian motion. This suggests to introduce the limiting process~$(\overline{X}_t)_{t \geq 0}$ obtained from~$(X_{\xi t})_{t \geq 0}$ as~$\xi \to +\infty$, and highlights the fact that the relaxation time is of order~$\xi$ in this regime. 

The above considerations suggest the following scalings for the convergence rate:
\begin{itemize}
\item for the Hamiltonian limit~$\xi \to 0$, some decay is observed over long times of order~$1/\xi$, which is consistent with the scaling~$\bar{\lambda} \sim \xi$ in Theorem~\ref{thm1new}, and also the scaling $(\xi t)^{-\sigma}$ for the denominator on the right hand side of the estimate in Theorem~\ref{thm2new};
\item for the overdamped limit~$\xi \to +\infty$, some decay is observed over long times of order~$\xi$, which is consistent with the scaling~$\bar{\lambda} \sim \xi^{-1}$ in Theorem~\ref{thm1new}. However, the situation is different for the estimates of Theorem~\ref{thm2new}, as the denominator on the right hand side of the estimate scales as~$(\xi^{-1-1/\sigma}t)^{-\sigma}$. Therefore, decay is observed over longer time scales of order~$\xi^{1+1/\sigma}$. This is due to the fact that the collision term is weaker. Of course, the timescale comes closer to the usual time scale of order~$\xi$ as~$\sigma$ increases (recall that this parameter is determined by the condition that $\mathscr{G}^\sigma \in \mathrm{L}^1(\gamma)$, see~\eqref{eq:weightedpoi}). 
\end{itemize}
The scaling $\max(\xi,\xi^{-1})$ for the rate $\lambda$ in Theorem~\ref{thm1new} is sharp, and it is attained, e.g.~, for \eqref{VOU} with~$\mathscr{X} = \mathbb T^d$, $\phi\equiv0$ and $\psi(v) = |v|^2/2$. On the other hand, the optimality of the dependence from $\xi$ in the convergence estimate of Theorem~\ref{thm2new} is unclear.  
%---------------------------------------
\section{Notation and preliminary technical results}
\label{sec:preliminary}

Recall that we denote by $A^\star$ the adjoint of a closed operator $A$ on the Hilbert space~$\mathrm{L}^2(\Theta)$, which, for any smooth functions~$h_1,h_2$ with compact supports, satisfies
\begin{equation}
  \label{VOUo}
  \int_{\mathcal{X} \times \mathcal{V}} (Ah_1) h_2 \, d\Theta = \int_{\mathcal{X} \times \mathcal{V}} h_1 \left(A^\star h_2\right) d\Theta.
\end{equation}
The action of~$A^\star$ is obtained by integration by parts for differential operators. In particular, for $1 \leq i \leq d$,
\[
\partial_{x_i}^\star = -\partial_{x_i} + \partial_{x_i} \phi,
  \qquad
  \partial_{v_i}^\star = -\partial_{v_i} + \partial_{v_i} \psi.
\]
This allows to rewrite the operators appearing in~\eqref{VOU} as the sum of a skew-symmetric operator
\[
T = \nabla_v \psi \cdot \nabla_x - \nabla_x \phi \cdot \nabla_v = \sum_{i=1}^d \partial_{v_i}^\star\partial_{x_i} - \partial_{x_i}^\star\partial_{v_i}, 
\]
and a symmetric one
\begin{equation}
  \label{eq:S_decomposition}
  \mathcal{S} = -\nabla_v^\star \nabla_v = -\sum_{i=1}^d \partial_{v_i}^\star\partial_{v_i},
\end{equation}
so that~\eqref{VOU} rewrites
\[
\partial_t h + Th = \xi \mathcal{S} h.
\]
Under Assumption~\ref{ass3}, the operator~$-T + \xi \mathcal{S}$ is hypoelliptic (see~\cite{hormander1967hypoelliptic,villani2009hypocoercivity}), hence the solutions to~\eqref{VOU} are smooth. 

Solutions~$h$ to \eqref{VOU} live in the Bochner space~$\mathrm{L}^\infty(0,\infty;\mathrm{L}^2(\Theta))$ and are such that $\nabla_v h \in \mathrm{L}^2(0,\infty;\mathrm{L}^2(\Theta))$, as we show thanks to the following decay estimate:
\begin{equation}
  \label{energy}
  \frac{d}{dt} \left( \|h(t)\|^2_{\mathrm{L}^2(\Theta)} \right) = 2\left\langle h(t), Th(t)\right\rangle_{\mathrm{L}^2(\Theta)} + 2\xi \left\langle h(t), \mathcal{S}h(t)\right\rangle_{\mathrm{L}^2(\Theta)}= - 2\xi \| \nabla_v h(t)\|^2_{\mathrm{L}^2(\Theta)} \leq 0,
\end{equation}
in view of the skew-symmetry of~$T$ and~\eqref{eq:S_decomposition}. Then, if $h_0 \in \mathrm{L}^2(\theta)$, \eqref{energy} implies that~$h \in \mathrm{L}^\infty(0,\infty;\mathrm{L}^2(\Theta))$. Moreover, for any~$\tau \in \mathbb{R}_+$,
\[
2 \xi \int_0^\tau \|\nabla_v h(t)\|^2_{\mathrm{L}^2(\Theta)} \, dt = \|h(0)\|^2_{\mathrm{L}^2(\Theta)} - \|h(\tau)\|^2_{\mathrm{L}^2(\Theta)} \leq \|h(0)\|^2_{\mathrm{L}^2(\Theta)},
\]
hence, by letting~$\tau \to \infty$, it follows that~$\nabla_v h \in \mathrm{L}^2(0,\infty; \mathrm{L}^2(\Theta))$. 

The first result provides estimates on the skew symmetric part of the space-time Fokker--Planck operator~$\partial_t + T$ in terms of its symmetric part~$\mathcal{S}$, for solutions of~\eqref{VOU}.

\begin{lemma}
  Consider a solution of~\eqref{VOU}. Then, 
  \begin{equation}\label{regest}
    \|(\partial_t + T)h \|_{\mathrm{L}^2(\mathrm{U}_\tau \otimes \mu ; \mathrm{H}^{-1}(\gamma))} \leq \xi \| \nabla_v h\|_{\mathrm{L}^2(\mathrm{U}_\tau \otimes \Theta)}.
  \end{equation}
  In particular, solutions to~\eqref{VOU} belong to $\mathrm{H}^1_{\mathrm{kin}}.$
\end{lemma}
Note that the right-hand-side of~\eqref{regest} is finite since~$\nabla_v h \in \mathrm{L}^2(0,\infty; \mathrm{L}^2(\Theta))$.
\begin{proof}
Fix a smooth and compactly supported test function $u$ with $\|u\|_{\mathrm{L}^2(\mathrm{U}_\tau \otimes \mu ; \mathrm{H}^1(\gamma))} \leq 1$. Then, denoting by~$\langle \, \cdot,\cdot \,\rangle$ the duality bracket in the sense of distributions, we have
\[
\begin{aligned}
  \left\langle (\partial_t  + T)h , u \right\rangle & = \xi \left\langle \mathcal{S} h , u \right\rangle = - \xi \int_0^\tau \int_{\mathcal{X} \times \mathcal{V}} \nabla_v h \cdot \nabla_v u \, d\mathrm{U}_\tau \, d\Theta \\
  & \leq \xi \| \nabla_v h\|_{\mathrm{L}^2(\mathrm{U}_\tau \otimes \Theta)} \| \nabla_v u\|_{\mathrm{L}^2(\mathrm{U}_\tau \otimes \Theta)} \\
  & \leq \xi \| \nabla_v h\|_{\mathrm{L}^2(\mathrm{U}_\tau \otimes \Theta)} \| u\|_{\mathrm{L}^2(\mathrm{U}_\tau \otimes \mu,\mathrm{H}^1(\gamma))}.
  \end{aligned}
\]
The desired estimate follows by extending the estimate to functions in~$\mathrm{L}^2(\mathrm{U}_\tau \otimes \mu ; \mathrm{H}^1(\gamma))$ by density, and taking the supremum over this functional space.
\end{proof}

The next lemma, whose proof is omitted, provides estimates similar to \cite[Lemma~2.3]{cao2023explicit}, obtained with an adaptation of the approach of \cite[Lemmas~3.6 and~3.7]{bernard2020hypocoercivity} to space-time operators. It provides a fully explicit dependence of the estimates on the dimension~$d.$

\begin{lemma}
  \label{lemmbfls}
  Assume that $\phi$ is a smooth potential such that~$\mathrm{e}^{-\phi} \in \mathrm{L}^1(\mathcal{X})$ with $\int_\mathcal{X} \mathrm{e}^{-\phi(x)} \,dx = 1$, and that~\eqref{poiphi}-\eqref{hess} hold true. Then, Assumption \ref{ass1bis} holds true with 
  \[
  L_\phi = 2\max \left(2 , \sqrt{d c''_\phi} \right).
  \]
  Moreover, for any $z \in \mathrm{H}^2(\mathrm{U}_\tau \otimes \mu)$ such that $\nabla_{t,x} z \in \mathrm{H}^1_{\mathrm{DC}}(\mathrm{U}_\tau \otimes \mu)^d$,
  \begin{equation}
    \label{bochnerx}
    \left\|\nabla_{t,x}^2 z \right\|_{\mathrm{L}^2\left(\mathrm{U}_\tau \otimes \mu\right)}^2 \leq 2 \|\nabla_{t,x}^\star \nabla_{t,x} z \|^2_{\mathrm{L}^2(\mathrm{U}_\tau \otimes \mu)} + 2 c'_\phi \left(\sqrt{d} + 2 \max \left(8 c'_\phi , \sqrt{c''_\phi \,d} \right)  \,  \right) \,  \|\nabla_x z \|^2_{\mathrm{L}^2\left(\mathrm{U}_\tau \otimes \mu\right)}.
\end{equation}
\end{lemma}

Note that the last term involves only derivatives in~$x$ and not in~$t$. It would be possible to have better constants in \eqref{bochnerx} under additional assumptions, for example if $\phi$ is strongly convex or has a Hessian bounded from below, see~\cite{bernard2020hypocoercivity,cao2023explicit}. 

%---------------------------------------
\section{Modified Poincar\'e inequalities and exponential decay rates}\label{sec:sec3}

We start by proving Lemma~\ref{avglem} in Section~\ref{sec:avglem}. We next turn to the proof of Theorem~\ref{thm1new} in Section~\ref{sec:proof_thm1}, for which a key element is a space-time Poincar\'e-type inequality (see Proposition~\ref{prop:Poincare_type}). We finally prove Corollary \ref{cor1new} in Section~\ref{subseccor1}. 

%---------------------------
\subsection{Proof of Lemma \ref{avglem}}\label{sec:avglem}

We write in this section the proof of Lemma~\ref{avglem}. Some integrals in the calculation below are formal and represent duality products. In order to estimate the~$\mathrm{H}^{-1}$ norm of~$\nabla_{t,x} \Pi h$, we first compute the $\mathrm{H}^{-1}$ norm of~$\partial_t \Pi h$, and then turn to~$\nabla_x \Pi h$. To bound~$\|\partial_t \Pi h\|_{\mathrm{H}^{-1}(\mathrm{U}_\tau \otimes \mu)}$, we fix a test function~$z \in \mathrm{H}^1_\ZDCT(\mathrm{U}_\tau \otimes \mu)$ with $\|z\|^2_{\mathrm{H}^1(\mathrm{U}_\tau \otimes \mu)} \leq 1$; while the estimation of~$\|\nabla_x \Pi h\|_{\mathrm{H}^{-1}(\mathrm{U}_\tau \otimes \mu)}$ can be performed by taking a test function~$Z =(Z_1,\dots,Z_d)\in \mathrm{H}^1_\ZDCT(\mathrm{U}_\tau \otimes \mu)^d$ with $\|Z\|^2_{\mathrm{H}^1(\mathrm{U}_\tau \otimes \mu)} \leq 1$.

\subsubsection*{Estimate of the time derivative}
  For~$\partial_t \Pi h$, we start by noticing that $\nabla_v \Pi h = 0$ and
  \[
  \int_{\mathcal{V}} \nabla \psi(v) \, \gamma(dv) = 0,
  \]
  so that, using the skew-symmetry of $\partial_t + T$ and the Dirichlet boundary conditions in the time variable for~$z$ in the integration by parts for the second term, 
  \[
  \begin{aligned}
    & \int_0^\tau \int_{\mathcal{X}} (\partial_t \Pi h) \, z \, d\mathrm{U}_\tau \, d\mu = \int_0^\tau \int_{\mathcal{X}} \int_{\mathcal{V}} ((\partial_t + T) \Pi \, h)(t,x,v) \, z(t,x) \, \mathrm{U}_\tau(dt)\, \mu(dx) \,\gamma(dv) \\
    & \qquad = \int_0^\tau \int_{\mathcal{X}} \int_{\mathcal{V}} [(\partial_t + T) h]  \, z \, d\mathrm{U}_\tau \, d\Theta + \int_0^\tau \int_{\mathcal{X}} \int_{\mathcal{V}} [(\Id-\Pi)h] \, (\partial_t + T) z \, d\mathrm{U}_\tau \, d\Theta \\
    & \qquad \leq \|(\partial_t + T) h\|_{\mathrm{L}^2(\mathrm{U}_\tau \otimes \mu,\mathrm{H}^{-1}(\gamma))} \|z\|_{\mathrm{L}^2(\mathrm{U}_\tau \otimes \mu)} \|1\|_{\mathrm{H}^1(\gamma)} + \|(\Id-\Pi)h\|_{\mathrm{\mathrm{L}^2}(\mathrm{U}_\tau \otimes \Theta)} \|(\partial_t + T)z\|_{\mathrm{L}^2(\mathrm{U}_\tau \otimes \Theta)}.
  \end{aligned}
  \]
  Denoting by~$\rho(\mathcal{M)} > 0$ the largest eigenvalue of the symmetric definite positive matrix~$\mathscr{M}$, the last factor in the last inequality can be bounded as
  \[
  \|(\partial_t + T)z\|_{\mathrm{L}^2(\mathrm{U}_\tau \otimes \Theta)} \leq \| \partial_t z \|_{\mathrm{L}^2(\mathrm{U}_\tau \otimes \mu)} + \| T z\|_{\mathrm{L}^2(\mathrm{U}_\tau \otimes \Theta)} \leq 1 + \sqrt{\rho(\mathscr{M})},
  \]
  where we used $\|z\|_{\mathrm{H}^1(\mathrm{U}_\tau \otimes \mu)} \leq 1$, and (recalling the definition~\eqref{eq:matrix_M} of~$\mathscr{M}$)
  \[
  \| T z\|_{\mathrm{L}^2(\mathrm{U}_\tau \otimes \Theta)}^2 = \| \nabla \psi \cdot \nabla_x z\|_{\mathrm{L}^2(\mathrm{U}_\tau \otimes \Theta)}^2 = \int_0^\tau \int_{\mathcal{X}} \nabla_x z \, \cdot \,  \mathscr{M} \, \nabla_x z \, d\mathrm{U}_\tau \, d\mu \leq \rho(\mathscr{M}) \|\nabla_x z\|_{\mathrm{L}^2(\mathrm{U}_\tau \otimes \mu)}^2.
  \]
  By taking the supremum over functions~$z \in \mathrm{H}^1_\ZDCT(\mathrm{U}_\tau \otimes \mu)$ with~$\|z\|^2_{\mathrm{H}^1(\mathrm{U}_\tau \otimes \mu)} \leq 1$, we finally obtain
  \begin{equation}
    \label{eq:dt_Pi_h_H-1}
    \|\partial_t \Pi h\|_{\mathrm{H}^{-1}(\mathrm{U}_\tau \otimes \mu)} \leq \|(\partial_t + T) h\|_{\mathrm{L}^2(\mathrm{U}_\tau \otimes \mu,\mathrm{H}^{-1}(\gamma))} + \left(1 + \sqrt{\rho(\mathscr{M})}\right)\|(\Id-\Pi)h\|_{\mathrm{\mathrm{L}^2}(\mathrm{U}_\tau \otimes \Theta)} .
  \end{equation}

  \subsubsection*{Estimate of the spatial derivatives}
  We next turn to the estimation of~$\|\nabla_x \Pi h\|_{\mathrm{H}^{-1}(\mathrm{U}_\tau \otimes \mu)}$. We introduce
  \[
  G = \frac{\nabla_v \psi}{\|\nabla_v \psi\|_{\mathrm{L}^2(\gamma)}}.
  \] 
  In view of the definition~\eqref{eq:matrix_M}, it holds by construction that
  \begin{equation}
    \label{eq:matrix_M_2}
    \int_\mathcal{V} G(v) \otimes M G(v) \, \gamma(dv) = \Id_{d\times d}, \qquad M = \|\nabla_v \psi\|_{\mathrm{L}^2(\gamma)}^2 \mathscr{M}^{-1}. 
  \end{equation}
  We use $G$ to reintroduce the velocity variable and reconstruct the operator~$T$ in the computation below. This allows to bypass the algebraic manipulations performed in the proof of~\cite[Theorem~2]{cao2023explicit}, and makes explicit the choice of functions considered in~\cite[Lemma~3.1]{albritton2024variational}. 

  The quantity to estimate is 
  \begin{equation}
    \label{eq:series_eq_nabla_x_Pi_Z}
  \begin{aligned}
   & \int_0^\tau \int_{\mathcal{X}} \nabla_x \Pi h \cdot Z \, d\mathrm{U}_\tau \, d\mu = \int_0^\tau \int_{\mathcal X} \int_{\mathcal{V}} \left[G(v) \cdot \nabla_x \Pi h(t,x)\right]\left[MG(v) \cdot Z(t,x)\right] \mathrm{U}_\tau(dt) \, \Theta(dx \, dv) \\
    & \qquad \qquad = \|\nabla_v \psi\|^{-1}_{\mathrm{L}^2(\gamma)}\int_0^\tau \int_{\mathcal X} \int_{\mathcal{V}} \left[(\partial_t + T)\Pi h\right](t,x,v)\left[MG(v) \cdot Z(t,x)\right] \mathrm{U}_\tau(dt) \, \Theta(dx \, dv) \\
    & \qquad \qquad = \|\nabla_v \psi\|^{-1}_{\mathrm{L}^2(\gamma)}\int_0^\tau \int_{\mathcal X} \int_{\mathcal{V}} \left[(\partial_t + T)h\right](t,x,v)\left[MG(v) \cdot Z(t,x)\right] \mathrm{U}_\tau(dt) \, \Theta(dx \, dv) \\
    & \qquad \qquad \ \ - \|\nabla_v \psi\|^{-1}_{\mathrm{L}^2(\gamma)}\int_0^\tau \int_{\mathcal X} \int_{\mathcal{V}} (T(1-\Pi) h)(t,x,v)\left[MG(v) \cdot Z(t,x)\right] \mathrm{U}_\tau(dt) \, \Theta(dx \, dv) \\
    & \qquad \qquad \ \ - \|\nabla_v \psi\|^{-1}_{\mathrm{L}^2(\gamma)}\int_0^\tau \int_{\mathcal X} \int_{\mathcal{V}} (\partial_t (1-\Pi) h)(t,x,v)\left[MG(v) \cdot Z(t,x)\right] \mathrm{U}_\tau(dt) \, \Theta(dx \, dv),
    \end{aligned}
  \end{equation}
  where the second equality follows from the fact that~$\nabla_v \Pi h = 0$ and $\int_{\mathcal{V}} G \, d\gamma = 0$. Since
  \[
  \|MZ\|_{\mathrm{L}^2(\mathrm{U}_\tau \otimes \mu)} \leq \rho(M) \|Z\|_{\mathrm{L}^2(\mathrm{U}_\tau \otimes \mu)} \leq \rho(M) \|Z\|_{\mathrm{H}^1(\mathrm{U}_\tau \otimes \mu)} \leq \rho(M),
  \]
  and~$\|G\|_{\mathrm{H}^1(\gamma)}$ is finite (in view of~\eqref{eq:integrability_conditions_on_psi}), the first integral in the last equality can be bounded by
  \[
  \|G \cdot MZ\|_{\mathrm{L}^2(\mathrm{U}_\tau \otimes \mu,\mathrm{H}^1(\gamma))} \|(\partial_t +T)h \|_{\mathrm{L}^2(\mathrm{U}_\tau \otimes \mu, \mathrm{H}^{-1}(\gamma) )} \leq \rho(M)\|G \|_{\mathrm{H}^1(\gamma)} \|(\partial_t +T)h \|_{\mathrm{L}^2(\mathrm{U}_\tau \otimes \mu, \mathrm{H}^{-1}(\gamma) )}.
  \]
  Since~$T$ is antisymmetric, the second integral in the last equality can be controlled by 
  \[
  \left| \int_0^\tau \int_{\mathcal X} \int_{\mathcal{V}} \left[(1-\Pi) h \right] \left[ T(MG \cdot Z)\right] d\mathrm{U}_\tau \, d\Theta \right| \leq \|(1-\Pi)h\|_{\mathrm{L}^2(\mathrm{U}_\tau \otimes \Theta)} \|T(MG \cdot Z)\|_{\mathrm{L}^2(\mathrm{U}_\tau \otimes \Theta)}.
  \]
  Introducing~$G^M = (G^M_1,\dots,G_d^M) = MG$ to simplify the notation,
  \[
  \begin{aligned}
    \left[ T (MG \cdot Z) \right](t,x,v)
    & = T\left[\sum_{i=1}^d G_i^M Z_i \right](t,x,v) \\
    & = \sum_{i=1}^d G_i^M(v) \nabla \psi(v) \cdot \nabla_x Z_i(t,x) - \nabla G_i^M(v) \cdot \nabla\phi(x) Z_i(t,x),
  \end{aligned}
  \]
  where~$G_i^M \nabla \psi \in \mathrm{L}^2(\gamma)^d$ and $\nabla G_i^M \in \mathrm{L}^2(\gamma)^d$ in view respectively of the first condition in~\eqref{eq:integrability_conditions_on_psi}; while~$Z_i \nabla \phi \in \mathrm{L}^2(\mathrm{U}_\tau \otimes \mu)$ with~$\|Z_i \nabla \phi \|_{\mathrm{L}^2(\mathrm{U}_\tau \otimes \mu)} \leq L_\phi \|Z_i \|_{\mathrm{H}^1(\mathrm{U}_\tau \otimes \mu)}$ by Assumption \ref{ass1bis}. The factor~$\|T(MG \cdot Z)\|_{\mathrm{L}^2(\mathrm{U}_\tau \otimes \Theta)}$ can therefore be controlled as
  \[
  \begin{aligned}
    \|T(MG \cdot Z)\|_{\mathrm{L}^2(\mathrm{U}_\tau \otimes \Theta)} & \leq \left( \sum_{i=1}^d \left\|G_i^M \nabla_v \psi\right\|^2_{\mathrm{L}^2(\gamma)} \right)^{1/2}\left(\sum_{i=1}^d \|\nabla_x Z_i\|_{\mathrm{L}^2(\mathrm{U}_\tau \otimes \mu)}^2\right)^{1/2} \\
    & \quad + L_\phi \left( \sum_{i=1}^d \left\|\nabla_v G_i^M \right\|^2_{\mathrm{L}^2(\gamma)}\right)^{1/2}\left(\sum_{i=1}^d \|Z_i\|_{\mathrm{L}^2(\mathrm{U}_\tau \otimes \mu)}^2\right)^{1/2} \\
    & \leq \left( \sum_{i=1}^d \left\|G_i^M \nabla_v \psi\right\|^2_{\mathrm{L}^2(\gamma)}\right)^{1/2} + L_\phi \left(\sum_{i=1}^d \left\|\nabla_v G_i^M \right\|^2_{\mathrm{L}^2(\gamma)}\right)^{1/2}.
    \end{aligned}
  \]
  The last integral in~\eqref{eq:series_eq_nabla_x_Pi_Z} can be bounded, after transferring the partial derivative in time to the test function, by~$\|(1-\Pi)h\|_{\mathrm{L}^2(\mathrm{U}_\tau \otimes \mu)} \|MG \cdot \partial_t Z\|_{\mathrm{L}^2(\mathrm{U}_\tau \otimes \Theta)}$. The second factor tensorises as~$\|MG\|_{\mathrm{L}^2(\gamma)} \| \partial_t Z\|_{\mathrm{L}^2(\mathrm{U}_\tau \otimes \mu)} \leq \rho(M) \|G\|_{\mathrm{L}^2(\gamma)} = \rho(M)$.

\subsubsection*{Conclusion of the proof}
By gathering the above estimates for the time and space derivatives, we finally obtain 
\begin{align*}
   &\left| \int_0^\tau \int_{\mathcal{X}} \nabla_{t,x} \Pi h \cdot (z,Z) \, d\mathrm{U}_\tau \, d\mu\right| \\
   &\leq \|(\partial_t + T) h\|_{\mathrm{L}^2(\mathrm{U}_\tau \otimes \mu,\mathrm{H}^{-1}(\gamma))} + \left(1 + \sqrt{\rho(\mathscr{M})} \right)\|(\Id-\Pi)h\|_{\mathrm{\mathrm{L}^2}(\mathrm{U}_\tau \otimes \Theta)} \\
   & \quad + \rho(M) \|\nabla_v \psi\|^{-1}_{\mathrm{L}^2(\gamma)} \|G \|_{\mathrm{H}^1(\gamma)} \|(\partial_t +T)h \|_{\mathrm{L}^2(\mathrm{U}_\tau \otimes \mu, \mathrm{H}^{-1}(\gamma) )} \\
   & \quad+ \|\nabla_v \psi\|^{-1}_{\mathrm{L}^2(\gamma)}  \|(1-\Pi)h\|_{\mathrm{L}^2(\mathrm{U}_\tau \otimes \Theta)} \left( \sqrt{\sum_{i=1}^d \left\|G_i^M \nabla_v \psi \right\|^2_{\mathrm{L}^2(\gamma)}} + L_\phi \sqrt{\sum_{i=1}^d \left\|\nabla G_i^M\right\|^2_{\mathrm{L}^2(\gamma)}}  \right) \\
  &\quad + \rho(M) \|\nabla_v \psi\|^{-1}_{\mathrm{L}^2(\gamma)} \|(1-\Pi)h\|_{\mathrm{L}^2(\mathrm{U}_\tau \otimes \Theta)}.
\end{align*}
This leads to the desired estimate, with the constant
\begin{equation}
  \label{kavg}
  K_\mathrm{avg} =  2\max\left(\mathscr C_1,\mathscr C_2\right)^2,
\end{equation}
where
\[
\begin{aligned}
  \mathscr C_1 & = 1 + \rho(M) \frac{\|G \|_{\mathrm{H}^1(\gamma)}}{\|\nabla_v \psi\|_{\mathrm{L}^2(\gamma)}}, \\
  \mathscr C_2 & = 1 + \sqrt{\rho(\mathscr{M})}  + \frac{\rho(M) + \sqrt{\displaystyle \sum_{i=1}^d \left\|G_i^M \nabla_v \psi \right\|^2_{\mathrm{L}^2(\gamma)}} + L_\phi \sqrt{\sum_{i=1}^d \left\|\nabla G_i^M \right\|^2_{\mathrm{L}^2(\gamma)}}}{\displaystyle\|\nabla_v \psi\|_{\mathrm{L}^2(\gamma)}}.
\end{aligned}
\]

\begin{remark}\label{rmk34}
Let us motivate here the choice of the function~$G$ in the previous calculation. It yields an orthogonal decomposition in~$\mathrm{L}^2(\mathrm{U}_\tau \otimes \Theta)$ (with respect to the scalar product induced by~$M$ in the~$v$ variable) of the term
\[
\left[(\partial_t + T) \Pi  h\right](t,x,v) =  (\partial_t \Pi h)(t,x) + G(v) \int_{\mathcal{V}} G(v') \cdot M(T \Pi h)(t,x,v') \gamma(dv').
\]
One component lies in the image of $\Pi$, while the second one is spanned by $G$ itself. Our strategy should be compared with the ones used to prove~\cite[Lemma 3.1]{albritton2024variational}, which is not quantitative, since the authors do not provide an explicit expression to terms like~$M$ and~$G$ in our proof. They can also be compared with the estimates at the beginning of \cite[Theorem 2]{cao2023explicit}. In spite of the similarity in the strategy, we make the algebra behind the proof more transparent and usable also for general kinetic energies $\psi$. 
\end{remark}

%--------------------------------
\subsection{Proof of Theorem~\ref{thm1new}}
\label{sec:proof_thm1}

The important Poincar\'e-type inequality which is used to prove Theorem~\ref{thm1new} is the following. Its proof crucially relies on Lemma \ref{avglem}. 

\begin{proposition}
  \label{prop:Poincare_type}
Under Assumptions~\ref{ass1}, \ref{ass2} and \ref{ass1bis}, it holds
\begin{equation}
  \label{modpoi}
  \forall h \in \mathrm{H}^1_{\mathrm{kin}} \cap \mathrm{L}^2_0(\mathrm{U}_\tau \otimes \Theta),
  \qquad 
  \bar{\lambda} \| h\|^2_{\mathrm{L}^2(\mathrm{U}_\tau \otimes\Theta)} \leq  \|(\partial_t + T)h\|^2_{\mathrm{L}^2(\mathrm{U}_\tau \otimes \mu;\mathrm{H}^{-1}(\gamma))} + \|(\Id- \Pi)h\|^2_{\mathrm{L}^2(\mathrm{U}_\tau \otimes \Theta)},
\end{equation}
with
\[
\bar{\lambda} = \frac{1}{1 + C_\tau^{\mathrm{Lions}} K_\mathrm{avg}},
\]
where $K_{\mathrm{avg}}$ is the constant appearing in Lemma \ref{avglem}.
\end{proposition}

\begin{proof}%[Proof of Proposition~\ref{prop:Poincare_type}]
  Since~$\Pi$ is an orthogonal projector, 
  \begin{equation}\label{wpoiint}
  \|h\|^2_{\mathrm{L}^2(\mathrm{U}_\tau \otimes \Theta)} = \|(\Id-\Pi)h\|^2_{\mathrm{L}^2(\mathrm{U}_\tau \otimes \mu)} + \|\Pi h\|^2_{\mathrm{L}^2(\mathrm{U}_\tau \otimes \mu)}.
  \end{equation}
  We apply~\eqref{lionstx} to the second term on the right hand side in \eqref{wpoiint}, and then use Lemma~\ref{avglem} to write (noting that~$\Pi h \in \mathrm{L}^2_0(\mathrm{U}_\tau \otimes \mu)$)
  \[
  \begin{aligned}
  \|\Pi h\|^2_{\mathrm{L}^2(\mathrm{U}_\tau \otimes \mu)} &\leq C^{\mathrm{Lions}}_\tau \|\nabla_{t,x} \Pi h\|^2_{\mathrm{H}^{-1}(\mathrm{U}_\tau \otimes \mu)} \\
  &\leq C^{\mathrm{Lions}}_\tau \mathrm{K}_{\mathrm{avg}} \left( \|(\Id-\Pi)h\|^2_{\mathrm{L}^2(\mathrm{U}_\tau \otimes \Theta)} + \|(\partial_t+T)h \|^2_{\mathrm{L}^2(\mathrm{U}_\tau \otimes \mu, \mathrm{H}^{-1}(\gamma))} \right).
  \end{aligned}
  \]
  Plugging the last estimate into \eqref{wpoiint} provides the claimed result.
\end{proof}

We are now in position to write the proof of Theorem~\ref{thm1new}. The exponential convergence is a direct consequence of \eqref{modpoi}. Indeed, note first that~\eqref{energy} implies
\begin{equation}\label{intpass}
\begin{aligned}
  \frac{d \mathscr{H}_\tau(t)}{dt}
  & =\frac{d}{dt} \left( \int_0^\tau \|  h(t+s,\cdot,\cdot)\|^2_{\mathrm{L}^2(\Theta)} \, ds \right) =  \int_0^\tau \frac{d}{dt} \left( \|  h(t+s,\cdot,\cdot)\|^2_{\mathrm{L}^2(\Theta)} \right) ds   \\
  & =-2\xi \int_0^\tau \| \nabla_v   h(t+s,\cdot,\cdot)\|^2_{\mathrm{L}^2(\Theta)} \, ds = -2 \xi \tau \int_0^\tau \| \nabla_v   h(t+s,\cdot,\cdot)\|^2_{\mathrm{L}^2(\Theta)} \, \mathrm{U}_\tau(ds).
\end{aligned}
\end{equation}
The space-time Poincar\'e inequality~\eqref{modpoi} implies, together with the bounds~\eqref{regest} and~\eqref{eq:coercivity_Delta_psi_OU}, that
\[
\bar{\lambda} \int_0^\tau \|    h(t+s,\cdot,\cdot)\|^2_{\mathrm{L}^2(\Theta)} \mathrm{U}_\tau(ds) \leq  \left(\frac{1}{c_\psi} + \xi^2\right) \int_0^\tau \| \nabla_v   h(t+s,\cdot,\cdot)\|^2_{\mathrm{L}^2(\Theta)} \mathrm{U}_\tau(ds).
\]
The combination of the previous inequality and \eqref{intpass} leads to
\[
\frac{d \mathscr{H}_\tau(t)}{dt} \leq -2  \bar{\lambda}\left(\frac{1}{\xi c_\psi} + \xi\right)^{-1}\mathscr{H}_\tau(t),
\]
from which  Theorem~\ref{thm1new} follows by a Gronwall inequality. 

%-----------------------
\subsection{Proof of Corollary \ref{cor1new}}\label{subseccor1}

Decay estimates on~$\mathscr{H}_\tau(t)$ lead to pointwise in time decay estimates on~$h(t)$ since~$t \mapsto \|h(t)\|^2_{\mathrm{L}^2(\Theta)}$ is nonincreasing, as a consequence of \eqref{energy}. Indeed, for any~$t \geq 0$,
\[
\| h(t+\tau) \|^2_{\mathrm{L}^2(\Theta)} \leq \frac1\tau \mathscr{H}_\tau(t) \leq \frac{\mathscr{H}_\tau(0)}{\tau} \rme^{-2\lambda t} \leq \| h(0) \|^2_{\mathrm{L}^2(\Theta)} \rme^{-2\lambda t},
\]
so that
\[
\forall t \geq \tau, \qquad \| h(t) \|_{\mathrm{L}^2(\Theta)} \leq \| h(0) \|_{\mathrm{L}^2(\Theta)} \rme^{-\lambda (t-\tau)}.
\]
The same bound holds for~$t \in [0,\tau]$ as~$\| h(t) \|_{\mathrm{L}^2(\Theta)} \leq \| h(0) \|_{\mathrm{L}^2(\Theta)}$. This finally gives the claimed result.

%-------------------------------------
\section{Algebraic decay rates and Lyapunov functionals}
\label{sec:sec4}

This section is devoted to the proof of Theorem~\ref{thm2new} and Corollary~\ref{cor2new}. The overall strategy is similar to the one presented in~\cite[Section~5]{brigati2022time}, so that we only sketch the main steps of the argument. Note however that, compared to~\cite{brigati2022time} (and various other works, for instance~\cite{bouin2019hypocoercivity}), the transport part of the evolution includes the terms~$\nabla \phi \cdot \nabla_v$ and~$\nabla \psi\cdot \nabla_x$ (as in~\cite{grothaus2019weak} for instance, instead of~$v \cdot \nabla_x$ as in~\cite{bouin2019hypocoercivity}).

\begin{proof}[Proof of Theorem~\ref{thm2new}]
The overall idea for the proof is to rely on the Poincar\'e-type estimate~\eqref{modpoi}, where the second term on the right hand side of this inequality is controlled by a combination of the weighted Poincar\'e inequality~\eqref{eq:weightedpoi} and a moment bound. More precisely, consider H\"older conjugate exponents~$p,q \geq 1,$ with $p = \sigma^{-1}(1+\sigma)>1$, so that $q = 1+\sigma.$ In particular, $q/p = \sigma$. 

\subsubsection*{Control of~$\|(\Id- \Pi)h\|^2_{\mathrm{L}^2(\mathrm{U}_\tau \otimes \Theta)}$}
We first claim that, for a smooth function~$g : \mathcal{X} \times \mathcal{V} \to \mathbb{R}$ with compact support,
\begin{equation}
  \label{eq:static_ineq_for_algebraic_decay}
\|g-\Pi g\|_{\mathrm{L}^2(\Theta)}^2 \leq P_\psi^{1/p} \|\nabla_v g \|_{\mathrm{L}^2(\Theta)}^{2/p} \left( \int_\mathcal{X} \int_{\mathcal{V}} \mathscr{G}(v)^\sigma |g(x,v)-\Pi g(x)|^2\,\Theta(dx\,dv) \right)^{1/q},
\end{equation}
recalling that $P_\psi$ is the Poincar\'e constant of \eqref{eq:weightedpoi}, in Assumption \ref{ass5}.
To prove this bound, we start from the following pointwise estimate, obtained with a H\"older inequality: for a given~$x \in \mathcal{X}$,
\[
\begin{aligned}
&\|g(x,\cdot)-\Pi g(x)\|_{\mathrm \mathrm{L}^2(\gamma)}^2 \\
&\qquad \leq \left(\int_{\mathcal{V}} \,\mathscr G(v)^{-1} |g(x,v)-\Pi g(x)|^2 \,\gamma(dv) \right)^{1/p} \left( \int_{\mathcal{V}} \mathscr G(v)^{\frac{q}{p}} |g(x,v)-\Pi g(x)|^2\,\gamma(dv) \right)^{1/q}.
\end{aligned}
\]
The inequality~\eqref{eq:static_ineq_for_algebraic_decay} then follows by bounding the first factor on the right hand side of the previous inequality with~\eqref{eq:weightedpoi}, integrating with respect to~$\mu(dx)$, and using a H\"older inequality in the~$x$ variable.

Consider next a smooth function~$h : \mathbb{R}_+ \times \mathcal{X} \times \mathcal{V} \to \mathbb{R}$ with compact support. Integrating the bound~\eqref{eq:static_ineq_for_algebraic_decay} for~$g(s) = h(t+s)$ with respect to~$\mathrm{U}_\tau(ds)$ and using a H\"older inequality in time then gives
\[
\begin{aligned}
& \int_0^\tau \|h(t+s,\cdot,\cdot)-\Pi \, h(t+s,\cdot)\|_{\mathrm \mathrm{L}^2(\Theta)}^2 \,\mathrm{U}_\tau(ds)\\
  &\qquad \qquad \leq P_\psi^{1/p} \,\left( \int_0^\tau \, \|\nabla_v h(t+s,\cdot,\cdot)\|_{\mathrm \mathrm{L}^2(\Theta)}^2\,\mathrm{U}_\tau(ds) \right)^{1/p} \\
  & \qquad \qquad \quad \times \left( \int_0^\tau \int_\mathcal{X} \int_\mathcal{V} \mathscr G(v)^{\sigma} \,|h(t+s,x,v)-\Pi h(t+s,x)|^2\, \Theta(dx \, dv) \, \mathrm{U}_\tau(ds) \right)^{1/q}.
 \end{aligned}
\]
When~$h$ solves \eqref{VOU} and~$h_0 \in \mathrm{L}^\infty(\mathcal{X} \times \mathcal{V})$, the strong maximum principle for hypoelliptic operators~\cite{taira2019strong} ensures that~$\|h(t)\|_{\mathrm{L}^\infty(\mathcal{X} \times \mathcal{V})} \leq \|h_0\|_{\mathrm{L}^\infty(\mathcal{X} \times \mathcal{V})}$, so that the integral in the last factor can be bounded by
\[
2\int_0^\tau \int_{\mathcal{X}} \int_{\mathcal{V}} \mathscr G^\sigma(v) \left( h^2(t+s,x,v) + (\Pi h)^2(t+s,x) \right) \, \Theta(dx\,dv) \, \mathrm{U}_\tau(ds) \leq 4 \left\| \mathscr{G}^\sigma \right\|_{\mathrm{L}^1(\gamma)} \|h_0\|^2_{\mathrm{L}^\infty(\mathcal{X} \times \mathcal{V})}.
\]
Finally, a density argument shows that
\[
\begin{aligned}
  & \int_0^\tau \|h(t+s,\cdot,\cdot)-\Pi \, h(t+s,\cdot)\|_{\mathrm \mathrm{L}^2(\Theta)}^2 \,\mathrm{U}_\tau(ds) \\
  & \qquad \qquad \leq P_\psi^{1/p} \left(4 \left\| \mathscr G^\sigma \right\|_{\mathrm{L}^1(\gamma)}\right)^{1/q} \|h_0\|^{2/q}_{\mathrm{L}^\infty(\mathcal{X} \times \mathcal{V})} \left( \int_0^\tau \|\nabla_v h(t+s,\cdot,\cdot)\|_{\mathrm \mathrm{L}^2(\Theta)}^2 \,\mathrm{U}_\tau(ds) \right)^{1/p}.
\end{aligned}
\]

\subsubsection*{Decay inequality on~$\mathscr{H}_\tau(t)$}
Combining the last bound with~\eqref{modpoi}, and then using~\eqref{regest} gives 
\[
\begin{aligned}
 &\bar{\lambda} \int_0^\tau \|h(t+s,\cdot,\cdot)\|^2_{\mathrm{L}^2( \Theta)} \, \mathrm{U}_\tau(ds) \\
 & \qquad \leq \int_0^\tau \|(\partial_t + T) \,h(t+s,\cdot,\cdot)\|^2_{\mathrm{L}^2(\mu;\mathrm{H}^{-1}(\gamma))} \, \mathrm{U}_\tau(ds)  \\
  &\qquad \quad + P_\psi^{1/p} \left(4 \left\| \mathscr G^\sigma \right\|_{\mathrm{L}^1(\gamma)}\right)^{1/q} \|h_0\|^{2/q}_{\mathrm{L}^\infty(\mathcal{X} \times \mathcal{V})} \left( \int_0^\tau \|\nabla_v h(t+s,\cdot,\cdot)\|_{\mathrm \mathrm{L}^2(\Theta)}^2 \,\mathrm{U}_\tau(ds) \right)^{1/p} \\
&\qquad \leq \xi^2 \, \int_0^\tau \|\nabla_v h(t+s,\cdot,\cdot)\|^2_{\mathrm{L}^2( \Theta)} \, \mathrm{U}_\tau(ds) \\ 
&\qquad \quad + P_\psi^{1/p} \left(4 \left\| \mathscr G^\sigma \right\|_{\mathrm{L}^1(\gamma)}\right)^{1/q} \|h_0\|^{2/q}_{\mathrm{L}^\infty(\mathcal{X} \times \mathcal{V})} \left( \int_0^\tau \|\nabla_v h(t+s,\cdot,\cdot)\|_{\mathrm \mathrm{L}^2(\Theta)}^2 \,\mathrm{U}_\tau(ds) \right)^{1/p}.
\end{aligned}
\]
In view of the definition \eqref{timeavg} and using~\eqref{intpass} to write
\[
\int_0^\tau \|\nabla_v h(t+s,\cdot,\cdot)\|_{\mathrm \mathrm{L}^2(\Theta)}^2 \,\mathrm{U}_\tau(ds) = -\frac{1}{2\xi \tau} \frac{d}{dt} \mathscr{H}_\tau(t),
\]
the latter inequality can be stated as 
\begin{equation}\label{diffeq}
\mathscr{H}_\tau(t) \leq \frac{\xi}{2 \, \bar{\lambda}} \,  \, \left(- \frac{d}{dt} \, \mathscr{H}_\tau(t)\right) +   \xi^{-\frac{1}{p}} \,  \, M_0 \, \left(-\frac{d}{dt} \, \mathscr{H}_\tau(t) \right)^{\frac{1}{p}},
\end{equation} 
where
\[
M_0 := 2^{(2-\sigma)/(1+\sigma)} \tau^{1/(\sigma+1)} \frac{1}{\bar{\lambda}} P_\psi^{1/p} \left\| \mathscr G^\sigma \right\|_{\mathrm{L}^1(\gamma)}^{1/q} \|h_0\|^{2/q}_{\mathrm{L}^\infty(\mathcal{X} \times \mathcal{V})}.
\]

The proof can now be concluded using classical arguments, relying on a generalization of the Gronwall inequality leading to algebraic decay rates. We however carefully keep track of the parameters influencing the convergence rate, by retaining the largest term in~\eqref{diffeq}. In order to make this precise, consider the function $$y \mapsto \vartheta(y) = \, \frac{\xi}{2\bar{\lambda}} \, y  + M_0 \left( \frac{y}{\xi} \right)^{\frac{1}{p}},$$
which is strictly increasing from $[0,\infty)$ to $[0,\infty)$, and hence invertible with a strictly increasing inverse.
Set $y_0 = \vartheta^{-1}(\mathscr{H}_\tau(0)).$
Note that 
$$\forall y \leq y_0, \qquad \vartheta(y) \leq \left( \xi^{-\frac{1}{p}} M_0 +  \, \frac{\xi}{2\bar{\lambda}} \, y_0^{\frac{p-1}{p}} \right) \, y^{\frac{1}{p}}.$$
Hence, as $t \mapsto \mathscr{H}_\tau(t)$ is nonincreasing and therefore~$\vartheta^{-1}\left(\mathscr{H}_\tau(t)\right) \leq y_0$, we can use the previous inequality to write
\[
\begin{aligned}
\mathscr{H}_\tau(t) & = \vartheta\left(\vartheta^{-1}\left(\mathscr{H}_\tau(t)\right)\right)\leq \left( \xi^{-\frac{1}{p}} M_0 +  \, \frac{\xi}{2\bar{\lambda}} \, y_0^{\frac{p-1}{p}} \right) \left(\vartheta^{-1}\left(\mathscr{H}_\tau(t)\right)\right)^{1/p} \\
& \leq \left( \xi^{-\frac{1}{p}} M_0 +  \,  \frac{\xi}{2\bar{\lambda}} \, y_0^{\frac{p-1}{p}} \right) \, \left( - \frac{d}{dt} \, \mathscr{H}_\tau(t)\right)^{\frac{1}{p}},
\end{aligned}
\]
where the last inequality follows from~\eqref{diffeq}. Applying the standard Bihari--LaSalle \cite{bihari1956generalization,lasalle1949uniqueness} argument, one proves, as in \cite{brigati2022time},
\begin{equation}\label{blsimplified}
\mathscr{H}_\tau(t) \leq \left( \mathscr{H}_\tau(0)^{1-p} + (p-1) \left( \xi^{-\frac{1}{p}} M_0 +  \frac{\xi}{2\bar{\lambda}} \, y_0^{\frac{p-1}{p}} \right)^{-p} t \right)^{-\frac{1}{p-1}}.
\end{equation}
Equation \eqref{blsimplified} can be rewritten as 
\[
\mathscr{H}_\tau(t) \leq \frac{\mathscr{H}_\tau(0)}{\left( 1 + \left( \xi^{-\frac{1}{p}} c_1 + \xi c_2 \right)^{-p} t \right)^{\frac{1}{p-1}}},
\]
with
\begin{equation}
  \label{eq:expressions_c1_c2}
  c_1 = \sigma^{\sigma/(1+\sigma)} M_0 \mathscr{H}_\tau(0)^{-1/(\sigma+1)}, \qquad c_2 = \sigma^{ \sigma/(1+\sigma) } \frac{\xi}{2\bar{\lambda}} y_0^{1/(1+\sigma)} \mathscr{H}_\tau(0)^{-1/(1+\sigma)}.
\end{equation}
The proof of Theorem~\ref{thm2new} is concluded by noting that~$p-1=\sigma^{-1}$.
\end{proof}

%\begin{remark}
%A more accurate estimate could be achieved by rewriting \eqref{diffeq} as
%$$\frac{d}{dt} \, \mathscr{H}_\tau(t) \leq - \zeta(\mathscr{H}_\tau(t)),$$
%with $\zeta = \vartheta^{-1}$ (the function $\vartheta$ being indeed invertible), and then integrating the latter inequality.
%\end{remark}

\begin{proof}[Proof of Corollary~\ref{cor2new}] The result is a direct consequence of Theorem~\ref{thm2new}, using the bounds
  \[
  \|h(t+\tau)\|_{\mathrm \mathrm{L}^2(\Theta)}^2 \leq \tau^{-1} \mathscr{H}_\tau(t), \qquad \tau^{-1} \mathscr{H}_\tau(0) \leq \|h_0\|_{\mathrm \mathrm{L}^2(\Theta)}^2,
  \]
  which follow from the fact that~$s \mapsto \|h(s)\|_{\mathrm \mathrm{L}^2(\Theta)}^2$ is non increasing.
\end{proof}

%--------------------------------------
\section{Weighted Poincar\'e--Lions inequalities}
\label{sec:sec5}

We present a constructive proof of Theorem~\ref{thm2}, based on a generalization of the manipulations performed in~\cite[Lemma~2.6]{cao2023explicit}. Compared to the results presented there, we extend the approach to the situation when the operator~$\nabla_x^\star \nabla_x$ has a spectral gap but not necessarily a discrete spectrum. We also present the argument in a more abstract manner, in order to highlight the key steps in the construction, and henceforth allow for extensions/adaptations to other situations.

\subsection{Reduction to solving a divergence equation}\label{ss:divred}

Recall that the space $\mathrm{H}^1_\DC(\mathrm{U}_\tau \otimes \mu)$ is defined in \eqref{dcdef}.
The Lions inequality~\eqref{lionstx} is implied by the following property (in fact, \cite{amrouche2015lemma} shows the equivalence between the two statements): there exists a positive constant $C_\tau^\mathrm{div}$ such that, for any~$f \in \mathrm{L}_0^2(\mathrm{U}_\tau \otimes \mu)$, there is a solution~$Z = (Z_{0},Z_{1},\dots,Z_d) \in \mathrm{H}^1_\DC(\mathrm{U}_\tau \otimes \mu)^{d+1}$ to
\begin{equation}
  \label{eq:div_eq_Z}
  -\partial_t Z_0 + \sum_{i=1}^d \partial_{x_i}^\star Z_i = f,
\end{equation}
which is such that 
\begin{equation}
  \label{elliptic}
  \|Z\|_{\mathrm{H}^1(\mathrm{U}_\tau \otimes \mu)} = \sqrt{\sum_{i=0}^d \|Z_i \|^2_{\mathrm{H}^1(\mathrm{U}_\tau \otimes \mu)}} \leq C_\tau^\mathrm{div} \|f\|_{\mathrm{L}^2(\mathrm{U}_\tau \otimes \mu)}. 
\end{equation}
Of course, there may be many solutions to the divergence-like equation~\eqref{eq:div_eq_Z}. When such a solution exists, we can rewrite the $\mathrm{L}^2$ norm of~$f$ using a duality bracket between~$\nabla_{t,x} f$ and~$Z$, relying crucially on the fact that~$Z \in \mathrm{H}^1_\DC(\mathrm{U}_\tau \otimes \mu)^{d+1}$, both for the integration by parts in the~$t$ variable, but also more generally to have a well defined duality product and transfer distributional derivatives in~$x$ from~$Z_i$ to~$f$:
%\todo[inline]{Just to confirm: If one was thinking only of integration by parts, it would seem that Dirichlet BCs in~$Z_0$ would be enough. Now, the way to think about the equalities below is that one should consider a $C^\infty$ function with compact support in~$(0,\tau) \times \mathcal{X}$, perform all manipulations in the sense of distributions, and pass to the~$H^1$ limit for~$f$. This is how Dirichlet boundary conditions arise for the other components of~$Z$. This should be maybe written this way?}
\[
\begin{aligned}
  \|f\|^2_{\mathrm{L}^2(\mathrm{U}_\tau \otimes \mu)}
  & = \int_0^\tau \int_\mathcal{X} f(t,x)^2 \, \mathrm{U}_\tau(dt) \, \mu(dx) = \int_0^\tau \int_\mathcal{X} \left(-\partial_t Z_0 + \sum_{i=1}^d \partial_{x_i}^\star Z_i\right) \, f \, d\mathrm{U}_\tau \, d\mu \\
  & = \left\langle \nabla_{t,x} f, Z\right\rangle_{\mathrm{H}^{-1}(\mathrm{U}_\tau \otimes \mu), H^{1}_\DC(\mathrm{U}_\tau \otimes \mu)} \leq \|Z\|_{\mathrm{H}^1(\mathrm{U}_\tau \otimes \mu)} \|\nabla_{t,x} f\|_{\mathrm{H}^{-1}(\mathrm{U}_\tau \otimes \mu)} \\
  & \leq C_\tau^\mathrm{div} \|f\|_{\mathrm{L}^2(\mathrm{U}_\tau \otimes \mu)} \|\nabla_{t,x} f\|_{\mathrm{H}^{-1}(\mathrm{U}_\tau \otimes \mu)}.
\end{aligned}
\]
Hence, Lions' inequality \eqref{lionstx} is recovered with~$C^{\mathrm{Lions}}_\tau \leq (C_\tau^\mathrm{div})^2$. We therefore concentrate in this section on obtaining solutions to~\eqref{eq:div_eq_Z} which satisfy~\eqref{elliptic}.   

%--------------------
\subsection{Solving the divergence equation}

We show in this section how to solve~\eqref{eq:div_eq_Z} in order to have estimates such as~\eqref{elliptic}. The precise result is the following.

\begin{proposition}
  \label{prop:estimates_divergence_equation}
  Suppose that \eqref{poiphi} and \eqref{hess} hold. Then, there exists an explicit constant~$C_\phi \in \mathbb{R}_+$, depending only on $c_\phi,c'_\phi,c''_\phi$ and not explicitly on $d$, such that, for any~$\tau > 0$ and~$f \in \mathrm{L}^2_0(\mathrm{U}_\tau \otimes \mu)$, the equation posed on the space-time domain~$[0,\tau] \times \mathcal{X}$
  \begin{equation}
    \label{eq:Lions_A}
    -\partial_t Z_0 + \sum_{i=1}^d \partial_{x_i}^\star Z_i = f,
  \end{equation}
  admits a solution~$(Z_0,\dots,Z_d) \in \mathrm{H}^1_{\DC}(\mathrm{U}_\tau \otimes \mu)^{d+1}$ which satisfies 
  \begin{equation}
\label{eq:estimates_solution_divergence_equation}
    \sqrt{\sum_{i=0}^d \| Z_i \|^2_{\mathrm{H}^1(\mathrm{U}_\tau \otimes \mu)}} \leq \, C_\tau^\mathrm{div} \, \|f\|_{\mathrm{L}^2(\mathrm{U}_\tau \otimes \mu)},
  \end{equation}
  for a constant $C_\tau^{\mathrm{div}}(\tau,d,C_\phi)>0$ independent from $f$.
  In addition, the following estimates hold
  \[
  C^{\mathrm{Lions}}_\tau \leq \left(C_\tau^\mathrm{div}\right)^2 \leq C_\phi \, \left( \frac{1}{\tau^2} + \left(1+\tau^2\right) \sqrt{d}\right). 
  \]
\end{proposition}

The explicit expression for $C_\phi$ can be deduced from~\eqref{eq:final_constant_Lions_ineq} in Section~\ref{sec:proof_estimates_divergence}. In fact, considering the positive self-adjoint operator satisfying (in view of \eqref{poiphi}) the following coercivity inequality on~$\mathrm{L}^2_0(\mu)$: 
\[
\Lop = \nabla_x^\star \nabla_x \geq c_\phi,
\]
and its square-root denoted by $L = \Lop^{1/2}$, we construct an explicit solution as follows by generalising the construction in \cite{cao2023explicit}: 
\begin{equation}
  \label{eq:explicit_solution_divergence_equation}
  Z = \nabla_{t,x} \mathscr{W}^{-1} \mathscr{P}_{\mathcal{N}^\perp} f + \begin{pmatrix} F_0(t,L) \\ \partial_{x_1} F_1(t,L) \\ \vdots \\ \partial_{x_d} F_1(t,L) \end{pmatrix} \mathscr{P}_{\mathcal{N},+} f + \begin{pmatrix} F_0(\tau-t,L) \\ \partial_{x_1} F_1(\tau-t,L) \\ \vdots \\ \partial_{x_d} F_1(\tau-t,L) \end{pmatrix} \mathscr{P}_{\mathcal{N},-} f,
\end{equation}
where we introduced the following objects:
\begin{itemize}
\item the operator~$\mathcal{W} = -\partial_t^2 + \Lop$ is considered on~$\mathrm{L}^2_0(\mathrm{U}_\tau \otimes \mu)$ and endowed with Neumann boundary conditions in time; 
\item the operator~$\mathscr{P}_{\mathcal{N}^\perp}$ is the projector onto the orthogonal of the vector space~$\mathcal{N} \subset \mathrm{L}^2_0(\mathrm{U}_\tau \otimes \mu)$ composed of linear combinations of forward and backward wave-like functions: $g \in \mathcal{N}$ if and only if there exist~$g_+,g_- \in \mathrm{L}^2_0(\mu)$ such that 
\begin{equation}
  \label{eq:decomposition_elements_N}
  g(t,x) = \left(\rme^{-t L}g_+\right)(x) + \left(\rme^{-(\tau-t) L} g_-\right)(x). 
\end{equation}
Note that~$\mathcal{W}g = 0$ in the weak sense by construction. Moreover, $g(t,\cdot) \in \mathrm{L}^2_0(\mu)$ for any~$t \in [0,\tau]$. 
\item the bounded operators~$F_0(t,L)$ and~$F_1(t,L)$ respectively act as
  \[
\begin{aligned}
  F_0(t,L) & = -2 L^{-1} \left(1-\rme^{-\tau L}\right)^{-2} \left(1-\rme^{-t L}\right)\left(1-\rme^{-(\tau-t) L}\right)\left(\rme^{-t L} - \frac{1+\rme^{-\tau L}}{2} \right), \\
  F_1(t,L) & = 6L^{-2} \left(1-\rme^{-\tau L}\right)^{-2} \left(1-\rme^{-t L}\right)\left(1-\rme^{-(\tau-t) L}\right) \, \rme^{-t L};
\end{aligned}
\]
\item the projectors~$\mathscr{P}_{\mathcal{N},+},\mathscr{P}_{\mathcal{N},-} : \mathrm{L}^2_0(\mathrm{U}_\tau \otimes \mu) \to \mathcal{N} \subset \mathrm{L}^2_0(\mathrm{U}_\tau \otimes \mu)$ act as
\begin{equation}
  \label{eq:def_Q_pm}
  \begin{aligned}
    \left(\mathscr{P}_{\mathcal{N},+} f\right)(t,x) & = \rme^{-t L} \left(\int_0^\tau G_+(s,L) \rme^{-s L} \, \mathrm{U}_\tau(ds) \right)^{-1} \int_0^\tau G_+(s,L)f(s,x) \, \mathrm{U}_\tau(ds), \\
    \left(\mathscr{P}_{\mathcal{N},-} f\right)(t,x) & = \rme^{-(\tau-t) L} \left(\int_0^\tau G_-(s,L) \rme^{-(\tau-s)L} \, \mathrm{U}_\tau(ds) \right)^{-1} \int_0^\tau G_-(s,L)f(s,x) \, \mathrm{U}_\tau(ds), 
\end{aligned}
\end{equation}
with
\begin{equation}
  \label{eq:U_tau_etc}
G_+(s,L) = V_\tau^{-1} \rme^{-(2\tau-s) L}-\rme^{-s L},
\qquad
G_-(s,L) = V_\tau^{-1} \rme^{-(\tau+s) L}-\rme^{-(\tau-s) L},
\qquad
V_\tau = \frac1\tau \int_0^\tau \rme^{-2 s L} \, ds.
\end{equation}
\end{itemize}
The form of the solution \eqref{eq:explicit_solution_divergence_equation} is motivated in Section~\ref{sec:heuristic_construction}, while the estimates~\eqref{eq:estimates_solution_divergence_equation} are proved in Section~\ref{sec:proof_estimates_divergence}. Note already that the operators~$\mathscr{P}_{\mathcal{N},+}$ and~$\mathscr{P}_{\mathcal{N},-}$ are bounded since~$\mathscr{P}_{\mathcal{N},+}^2 = \mathscr{P}_{\mathcal{N},+}$ and~$\mathscr{P}_{\mathcal{N},-}^2 = \mathscr{P}_{\mathcal{N},-}$. To obtain the estimates~\eqref{eq:estimates_solution_divergence_equation} for the last two terms in~\eqref{eq:explicit_solution_divergence_equation}, it is therefore sufficient to obtain uniform-in-time bounds on the operators~$F_0(t,L)$, $\partial_t F_0(t,L)$, $\partial_{x_i} F_0(t,L)$ and~$\partial_{x_i} F_1(t,L),\partial_t \partial_{x_i} F_1(t,L),\partial_{x_i,x_j}^2 F_1(t,L)$ considered on~$\mathrm{L}_0^2(\mu)$. 

An inspection of the proof finally shows that the factor~$\tau^2$ on the right hand side of~\eqref{eq:estimates_solution_divergence_equation} comes from the first term of the solution~\eqref{eq:explicit_solution_divergence_equation} (it can be traced back to some coercivity estimate based on the Poincaré inequality for~$\mathrm{U}_\tau$), while the factor~$\tau^{-2}$ arises from the part of the solution associated with~$\mathcal{P}_\mathcal{N,\pm}f$, because of the operator~$\left(1-\rme^{-\tau L}\right)^{-2}$ in the expressions of~$F_0(t,L)$ and~$F_1(t,L)$.

\subsection{Heuristic construction of the solution}
\label{sec:heuristic_construction}

To find one possible solution to~\eqref{eq:Lions_A}, a natural idea is to look for~$(Z_0,\dots,Z_d)$ of the form~$(\partial_t u, \partial_{x_1} u, \dots,\partial_{x_d} u)$, with~$u$ the unique solution of the following equation:
\begin{equation}
  \label{eq:eq_to_solve_with_LaxMilgram}
  (-\partial_t^2 + \Lop)u = f, \qquad \partial_t u(0,\cdot) = \partial_t u(\tau,\cdot) = 0.
\end{equation}
The functions~$\partial_{x_i} u(0,\cdot)$ and~$\partial_{x_i} u(\tau,\cdot)$ are however not~0 in general, so that $Z$ would not belong to $\mathrm{H}^1_\DC(\mathrm{U}_\tau \otimes \mu)$. The functions~$\partial_{x_i} u(0,\cdot)$ and~$\partial_{x_i} u(\tau,\cdot)$ turn out to vanish provided the right hand side of~\eqref{eq:Lions_A} belongs to the orthogonal of the vector space~$\mathcal{N}$.

The first step of the construction, as done in~\cite[Lemma~2.4]{cao2023explicit}, is to solve~\eqref{eq:Lions_A} for the part of the right hand side~$f$ which is in~$\mathcal{N}^\perp$, using~\eqref{eq:eq_to_solve_with_LaxMilgram} and the Lax--Milgram lemma (see~\cite[Chapter 5]{brezis2010functional} for instance); see Section \ref{sec:proof_estimates_divergence} for the proof of the following lemma.

\begin{lemma}
  \label{lem:Lax_Milgram}
  Suppose that~\eqref{poiphi} and~\eqref{hess} hold true.
  Fix~$\tau > 0$. For any~$f \in \mathrm{L}^2_0(\mathrm{U}_\tau \otimes \mu)$, there exists a unique weak solution~$u \in \mathrm{H}^1(\mathrm{U}_\tau \otimes \mu) \cap \mathrm{L}^2_0(\mathrm{U}_\tau \otimes \mu)$ to~\eqref{eq:eq_to_solve_with_LaxMilgram}: 
  \begin{equation}
    \label{eq:weak_formulation_LM}
    \forall w \in \mathrm{H}^1(\mathrm{U}_\tau \otimes \mu), \qquad 
  \int_0^\tau \int_\mathcal{X} \left[ (\partial_t u)(\partial_t w) + \nabla_x u \cdot \nabla_x w \right]\, d\mathrm{U}_\tau\, d\mu =\int_0^\tau \int_\mathcal{X} f w \, d\mathrm{U}_\tau \,d\mu.
  \end{equation}
  Moreover, this solution belongs to~$\mathrm{H}^2(\mathrm{U}_\tau \otimes \mu)$ and satisfies
  \[
  \|\nabla_{t,x} u \|^2_{\mathrm{H}^1(\mathrm{U}_\tau \otimes \mu)} \leq C_\tau^\mathrm{LM} \|f\|^2_{\mathrm{L}^2(\mathrm{U}_\tau \otimes \mu)},
  \]
  where
  \begin{equation}
    \label{eq:C_tau_LM}
  C_\tau^\mathrm{LM} = 2 + \frac{1 + 2 c'_\phi \left[\sqrt{d} + 2 \max \left( 8 c'_\phi , \sqrt{c''_\phi \,d} \right) \right]}{C_\mathrm{P}},
  \end{equation}
  with
  \begin{equation}
    \label{eq:def_C_P}
    C_{\mathrm P} = \min \left( c_\phi, \frac{\pi^2}{\tau^2}\right).
  \end{equation}
  Finally, when $f \in \mathcal{N}^\perp$, then $\nabla_{t,x} u \in \mathrm{H}^1_\DC(\mathrm{U}_\tau \otimes \mu)$.
\end{lemma}

The second step is to construct a solution to~\eqref{eq:Lions_A} for the part of the right hand side in~$\mathcal{N}$. This is done in two steps: (i) making explicit the expression of the orthogonal projection~$\mathscr{P}_\mathcal{N} = \mathscr{P}_{\mathcal{N},+} + \mathscr{P}_{\mathcal{N},-}$ of~$\mathrm{L}^2_0(\mathrm{U}_\tau \otimes \mu)$ onto~$\mathcal{N}$; and then (ii) solving the divergence equation~\eqref{eq:Lions_A} with right hand side~$\mathscr{P}_{\mathcal{N},\pm} f$. To identify~$\mathscr{P}_{\mathcal{N},\pm}$, we write 
\begin{equation}
  \label{eq:decomposition_P_N}
  (\mathscr{P}_\mathcal{N} f)(t,x) = \left(\rme^{-t L} Q_+ f\right)(x) + \left(\rme^{-(\tau-t) L} Q_- f\right)(x).
\end{equation}
By definition of the orthogonal projection, it holds, for any~$w_+,w_- \in \mathrm{L}^2_0(\mu)$, 
\begin{align*}
  \int_0^\tau \int_\mathcal{X} \left[f(t,x)-(\mathscr{P}_\mathcal{N} f)(t,x)\right]\left(\rme^{-t L} w_+\right)(x) \, \mathrm{U}_\tau(dt) \, \mu(dx) & = 0, \\
  \int_0^\tau \int_\mathcal{X} \left[f(t,x)-(\mathscr{P}_\mathcal{N} f)(t,x)\right]\left(\rme^{-(\tau-t) L} w_-\right)(x) \, \mathrm{U}_\tau(dt) \, \mu(dx) & = 0.
\end{align*}
Therefore, using the self-adjointness of $L,$
\begin{equation}
  \label{eq:system_for_projection}
  \begin{aligned}
    &\int_\mathcal{X} (V_\tau Q_+ f) w_+ \, d\mu + \int_\mathcal{X} \left(\rme^{-\tau L} Q_- f\right) w_+ \, d\mu  = \int_\mathcal{X} \left[\int_0^\tau \rme^{-t L} f(t,x) \, \mathrm{U}_\tau(dt) \right] w_+(x) \, \mu(dx), \\
    &\int_\mathcal{X} \left(\rme^{-\tau L} Q_+ f\right) w_- \, d\mu + \int_\mathcal{X} (V_\tau Q_- f) w_- \, d\mu = \int_\mathcal{X} \left[\int_0^\tau \rme^{-(\tau-t) L} f(t,x) \, \mathrm{U}_\tau(dt) \right] w_-(x) \, \mu(dx), \\ 
  \end{aligned}
\end{equation}
where we recall from~\eqref{eq:U_tau_etc} that
\[
V_\tau = \int_0^\tau \rme^{-2t L} \, \mathrm{U}_\tau(dt) = \frac{1-\rme^{-2\tau L}}{2\tau} L^{-1}.
\]
Since~$L \geq \sqrt{c_\phi}$ on~$\mathrm{L}^2_0(\mu)$ (in particular, $0 \leq \mathrm{e}^{-2\tau L} \leq \mathrm{e}^{-2 \tau \sqrt{c_\phi}}$ in the sense of symmetric operators and hence $1-\mathrm{e}^{-2\tau L}$ is invertible), the operator~$2\tau(1-\rme^{-2\tau L})^{-1}L$ is coercive on~$\mathrm{L}^2_0(\mu)$, hence invertible on this functional space, with inverse~$V_\tau$. Choosing~$w_+ = V_\tau^{-1} \rme^{-\tau L} w_-$ with $w_- \in \mathrm{L}^2_0(\mu) \cap D(L)$, and subtracting the second equation of~\eqref{eq:system_for_projection} from the first one, leads to
\[
\int_\mathcal{X} \left\{ \left[\left(V_\tau^{-1}\rme^{-2\tau L} - V_\tau\right) Q_- f\right] - \left[\int_0^\tau \left(V_\tau^{-1} \rme^{-(\tau+t) L}-\rme^{-(\tau-t)L}\right) \, f \, d\mathrm{U}_\tau \right] \right\} \, w_- \, d\mu =0 .
\]
Since the latter equality holds for any $w_- \in \mathrm{L}^2_0(\mu) \cap D(L)$, it follows that   
\begin{equation}
  \label{eq:formula_Q_minus}
  Q_- f(x) = \left(V_\tau^{-1}\rme^{-2\tau L} - V_\tau\right)^{-1} \int_0^\tau \left(V_\tau^{-1} \rme^{-(\tau+t) L}-\rme^{-(\tau-t) L}\right) f(t,x) \, \mathrm{U}_\tau(dt),
\end{equation}
which leads to~\eqref{eq:def_Q_pm} in view of~\eqref{eq:decomposition_P_N}. The formula for~$Q_+ f$ is obtained in a similar manner, upon replacing~$\tau-s$ by~$s$ everywhere. In order to prove that the operators $Q_{\pm}$ are well-defined, it remains to show that the operator in front of the time integral in~\eqref{eq:formula_Q_minus} is invertible on~$\mathrm{L}^2_0(\mu)$. This follows from the fact that this operator can be rewritten, up to a minus sign, as the inverse of
\[
V_\tau - V_\tau^{-1}\rme^{-2\tau L} = \frac{1}{2\tau} L^{-1} \left(1-\mathrm{e}^{-2\tau L}\right)^{-1} B(2\tau L),
\qquad
B(z) = \left(1-\mathrm{e}^{-z}\right)^2 - z^2 \mathrm{e}^{-z},
\]
with $B(2\tau L)$ coercive on~$\mathrm{L}^2_0(\mu)$. Indeed,
\[
B'(z) = 2\rme^{-z} \left( 1-z+\frac{z^2}{2}-\mathrm{e}^{-z}\right) \geq 0,
\]
so that $B(2\tau L) \geq B(2\tau \sqrt{c_\phi})$ with~$B(2\tau \sqrt{c_\phi})>0$.

Let us next solve  equation~\eqref{eq:Lions_A} with right hand side~$(\mathscr{P}_{\mathcal{N},+}f)(t) = \rme^{-t L} f_+$ for~$f_+ \in \mathrm{L}^2_0(\mu)$. Note that $(\mathscr{P}_{\mathcal{N},+}f)(t) \in \mathrm{L}^2_0(\mu)$ for all~$t \in [0,\tau]$. The manipulations performed here are similar for the right hand side~$\rme^{-(\tau-t) L} f_-$ with~$f_- \in \mathrm{L}^2_0(\mu)$, upon replacing~$t$ by~$\tau-t$ in the operators which appear. We consider the following ansatz, using functional calculus on the self-adjoint operator~$L = \Lop^{1/2}$, for functions~$\cF_{0},\cF_{1}:[0,\tau] \times [c_\phi,+\infty[ \to \mathbb{R}$ to be determined:
\[
Z_{0,+}(t) = \cF_{0}(t,L)f_+, \qquad Z_{i,+}(t) = \partial_{x_i} \cF_{1}(t,L) f_+.
\]
Note that these functions satisfy the equation~\eqref{eq:Lions_A} with right hand side~$\rme^{-t L}f_+$ provided the following equation is satisfied for operators on~$\mathrm{L}^2_0(\mu)$:
\[
-\partial_t \cF_{0}(t,L) + L^2 \cF_{1}(t,L) = \rme^{-t L}.
\]
Moreover, in view of establishing the estimates~\eqref{eq:estimates_solution_divergence_equation} and in order to make use of the fact that~$\mathscr{P}_{\mathcal{N},+}$ is bounded, it is convenient to rewrite~$\cF_{0}(t,L)f_+$ and~$\cF_1(t,L)f_+$ for~$t \in [0,\tau]$ given as the operators on~$\mathrm{L}^2_0(\mu)$
\begin{equation}
  \label{eq:from_F_to_mathcal_F}
  F_0(t,L) = \mathcal{F}_0(t,L) \rme^{t L}, \qquad F_1(t,L) = \mathcal{F}_1(t,L) \rme^{t L},
\end{equation}
applied to the function~$(\mathscr{P}_{\mathcal{N},+} f)(t) = \rme^{-t L}f_+$. This motivates looking for
\begin{equation}
  \label{eq:mathcal_F_as_polynomial}
  \cF_0(t,L) = P_{0}\left(\rme^{-t L}\right), \qquad \cF_1(t,L) = P_1\left(\rme^{-t L}\right),
\end{equation}
where~$P_0,P_1$ are polynomial functions satisfying the operator equation
\begin{equation}
  \label{eq:eq_P0_P1}
  -\frac{d}{dt} \left[ P_0\left(\rme^{-t L}\right) \right] + L^2 P_1\left(\rme^{-t L}\right) = \rme^{-t L},
\end{equation}
together with the boundary conditions
\begin{equation}
  \label{eq:BC_conditions_P0_P1}
  P_0(1) = P_0\left(\rme^{-\tau L}\right) = P_1(1) = P_1\left(\rme^{-\tau L}\right) = 0.
\end{equation}
In view of the boundary conditions~\eqref{eq:BC_conditions_P0_P1}, and the presence of the operators~$L$ and~$L^2$ in~\eqref{eq:eq_P0_P1}, a natural ansatz is
\[
\left\{ \begin{aligned}
  P_0\left(\rme^{-t L}\right) & = A_1 L^{-1} \left(1-\rme^{-t L}\right)\left(\rme^{-t L}-\rme^{-\tau L}\right)\left(\rme^{-t L} + A_2\right)\widetilde{P}_0\left(\rme^{-t L}\right), \\
  P_1\left(\rme^{-t L}\right) & = A_3 L^{-2} \left(1-\rme^{-t L}\right)\left(\rme^{-t L}-\rme^{-\tau L}\right)\widetilde{P}_1\left(\rme^{-t L}\right),
\end{aligned} \right.
\]
for operators~$A_1,A_2,A_3$ and polynomials~$\widetilde{P}_0,\widetilde{P}_1$ to be determined. This ansatz makes sense for operators considered on~$\mathrm{L}^2_0(\mu)$, as~$L$ is invertible on that functional space.

Note already that the operators~$P_0\left(\rme^{-t L}\right)$ and~$P_1\left(\rme^{-t L}\right)$ contain a composition with~$\rme^{-t L}-\rme^{-\tau L}$, which allows to further compose the operators with~$\rme^{t L}$ in order to make sense of~\eqref{eq:from_F_to_mathcal_F}. In order to find the operators~$A_1,A_2,A_3$, which are assumed to be functions of~$L$ (and hence commute with all operators at hand), we evaluate~\eqref{eq:eq_P0_P1}: denoting by~$\theta_\star = \rme^{-\tau L}$ and choosing~$\widetilde{P}_0 = 1$,
\[
A_1 \theta \frac{d}{d\theta}\left[ (1-\theta)(\theta-\theta_\star)(\theta+A_2)\right] + A_3 (1-\theta)(\theta-\theta_\star)\widetilde{P}_1(\theta) = \theta, 
\]
where the various functions are evaluated at $\theta = \rme^{-t L}$. This suggests to choose~$\widetilde{P}_1(\theta) = \theta$, so that
\[
A_1 \frac{d}{d\theta}\left[ (1-\theta)(\theta-\theta_\star)(\theta+A_2)\right] + A_3 (1-\theta)(\theta-\theta_\star) = \mathrm{Id}.
\]
The left hand side of this equation is a second order polynomial in~$\theta$, so that the operators~$A_1,A_2,A_3$ are obtained by identifying the prefactors of the various terms~$\theta^a$ for~$a = 0,1,2$. More precisely, % note \frac{d}{d\theta}\left[ (1-\theta)(\theta-\theta_\star)(\theta+A_2)\right] = -3\theta^2 + 2\theta(1+\theta_\star-A_2) + A_2(1+\theta_\star)-\theta_\star
\[
\left\{ \begin{aligned}
  A_1 \left[ A_2(1+\theta_\star) - \theta_\star \right] - A_3\theta_\star = 1, \\
  2A_1(1+\theta_\star-A_2) +A_3(1+\theta_\star) = 0, \\
  -3A_1-A_3 = 0.
\end{aligned} \right.
\]
The last condition implies that~$A_3 = -3A_1$. Plugging this equality in the second condition leads to $A_2 = -(1+\theta_\star)/2$. Finally, the first condition gives~$-A_1(1-\theta_\star)^2 = 1$, so that
\[
A_1 = -2\left(1-\rme^{-\tau L}\right)^{-2}.
\]
Combining~\eqref{eq:from_F_to_mathcal_F} and~\eqref{eq:mathcal_F_as_polynomial} then leads to the expressions of~$F_0,F_1$ in~\eqref{eq:explicit_solution_divergence_equation}.

%--------------
\subsection{Proof of the estimates~\eqref{eq:estimates_solution_divergence_equation}}
\label{sec:proof_estimates_divergence}

The estimate on the first term in~\eqref{eq:explicit_solution_divergence_equation} is a direct consequence of Lemma~\ref{lem:Lax_Milgram}.

\begin{proof}[Proof of Lemma~\ref{lem:Lax_Milgram}]
  The first part of the proof closely follows the one of~\cite[Lemma~2.4]{cao2023explicit}. We simply revisit it to make the constants precise. The existence and uniqueness of the weak solution~$u \in \mathrm{H}^1(\mathrm{U}_\tau \otimes \mu) \cap \mathrm{L}^2_0(\mathrm{U}_\tau \otimes \mu)$ to~\eqref{eq:weak_formulation_LM} is a consequence of Lax--Milgram's theorem on $\mathrm{H}^1(\mathrm{U}_\tau \otimes \mu) \cap \mathrm{L}_0^2(\mathrm{U}_\tau \otimes \mu)$ and the spatial Poincar\'e inequality~\eqref{poiphi}, which implies the following space-time Poincar\'e inequality by tensorisation:
  \[
  \forall v \in \mathrm{H}^1(\mathrm{U}_\tau \otimes \mu) \cap \mathrm{L}^2_0(\mathrm{U}_\tau \otimes \mu),
  \qquad
  \| v\|^2_{\mathrm{L}^2(\mathrm{U}_\tau \otimes \mu)} \leq \frac{1}{C_\mathrm{P}} \|\nabla_{t,x} v\|^2_{\mathrm{L}^2(\mathrm{U}_\tau \otimes \mu)},
  \]
  with~$C_\mathrm{P}$ defined in~\eqref{eq:def_C_P}. Moreover, 
  \begin{eqnarray}
    \label{eq:ellipt3}
    \|\nabla_{t,x} u\|^2_{\mathrm{L}^2(\mathrm{U}_\tau \otimes \mu)} \leq C_{\mathrm P}^{-1} \, \|f\|^2_{\mathrm{L}^2(\mathrm{U}_\tau \otimes \mu)}.
  \end{eqnarray}
  It can next be shown that~$u \in \mathrm{H}^2(\mathrm{U}_\tau \otimes \mu)$, with, thanks to~\eqref{bochnerx},
  \[
  \left\|\nabla^2_{t,x} u\right\|^2_{\mathrm{L}^2(\mathrm{U}_\tau \otimes \mu)} \leq 2 \|f\|^2_{\mathrm{L}^2(\mathrm{U}_\tau \otimes \mu)} + 2 c'_\phi \left[\sqrt{d} + 2 \max \left(8 c'_\phi , \sqrt{c''_\phi \,d} \right)  \right]  \|\nabla_{t,x} u\|^2_{\mathrm{L}^2(\mathrm{U}_\tau \otimes \mu)}.
  \]
   Hence, by \eqref{eq:ellipt3},
   \[
   \|\nabla_{t,x} u\|^2_{\mathrm{H}^1(\mathrm{U}_\tau\otimes\mu)} = \left\|\nabla^2_{t,x} u\right\|^2_{\mathrm{L}^2(\mathrm{U}_\tau \otimes \mu)} + \|\nabla_{t,x} u\|^2_{\mathrm{L}^2(\mathrm{U}_\tau \otimes \mu)} \leq C^\mathrm{LM}_\tau \, \|f\|^2_{\mathrm{L}^{2}(\mathrm{U}_\tau \otimes \mu)},
   \]
   with
   \[
   C^\mathrm{LM}_\tau =  2 + \frac{1 + 2 c'_\phi \left[\sqrt{d} + 2 \max \left(8 c'_\phi , \sqrt{c''_\phi \,d} \right) \right]}{C_\mathrm{P}}.
   \]

  We conclude by showing that the required boundary conditions are satisfied if $f \in \mathcal{N}^\perp$. 
  Consider an element~$w \in \mathcal{N}$, written as~\eqref{eq:decomposition_elements_N} for two smooth functions~$w_+,w_- \in \mathrm{L}^2_0(\mu)$ with compact support. Recall that~$(-\partial_t^2 + \Lop)w = 0$. Then, an integration by parts and the weak formulation~\eqref{eq:weak_formulation_LM} give
  \[
  \begin{aligned}
    \int_\mathcal{X} \left(u(\tau,\cdot)\partial_t w(\tau,\cdot) - u(0,\cdot)\partial_t w(0,\cdot)\right] d\mu & = \int_0^\tau \int_\mathcal{X} \left[ (\partial_t u)(\partial_t w) + \nabla_x u \cdot \nabla_x w \right]\, d\mathrm{U}_\tau \, d\mu \\
  & = \int_0^\tau \int_\mathcal{X} f w \, d\mathrm{U}_\tau \, d\mu = 0. 
  \end{aligned}
  \]
  the latter equality following from the fact that~$w \in \mathcal{N}$ and $f \in \mathcal{N}^\perp$ by assumption. Therefore, recalling~$L = \Lop^{1/2}$,  
  \[
  \int_\mathcal{X} \left[ u(\tau,\cdot) L\left(-\rme^{-\tau L}w_+ + w_-\right) - u(0,\cdot) L\left(-w_++\rme^{-\tau L}w_-\right)\right] d\mu = 0.
  \]
  Fixing~$w_-$ and choosing~$w_+ = \rme^{-\tau L}w_-$ leads to
  \[
  \int_\mathcal{X} u(\tau,\cdot) L\left(1-\rme^{-2\tau L}\right)w_- \, d\mu = 0.
  \]
  Since this equality holds true for any smooth~$w_- \in \mathrm{L}^2_0(\mu)$ with compact support, and $L \left(1-\rme^{-2\tau L}\right)$ is a coercive hence invertible operator on $\mathrm{L}^2_0(\mu)$, we can conclude that~$u(\tau,\cdot) = 0$. A similar reasoning allows to conclude that~$u(0,\cdot) = 0$. By taking spatial gradients (whose traces are well-posed, as $u \in \mathrm{H}^2(\mathrm{U}_\tau \otimes \mu)$), we finally obtain that~$\nabla_x u(0,\cdot) = \nabla_x u (\tau,\cdot) = 0$, and can therefore conclude that~$\nabla_{t,x} u \in \mathrm{H}^1_\DC(\mathrm{U}_\tau \otimes \mu)$ since the boundary conditions for the partial derivatives in time are satisfied by construction (recall~\eqref{eq:eq_to_solve_with_LaxMilgram}).
\end{proof}

We next show that the second and third terms in~\eqref{eq:explicit_solution_divergence_equation} satisfy~\eqref{eq:estimates_solution_divergence_equation}. Given that the two terms are extremely similar, we only treat the second one. We denote by~$\mathcal{B}(E)$ the Banach space of bounded linear operators on a Banach space~$E$.

\begin{lemma}\label{Nlem}
  Fix~$f \in \mathrm{L}^2_0(\mu)$. Under the same assumptions as in Theorem~\ref{thm2}, the function \[Z^+ = \left( F_0(t,L), \partial_{x_1} F_1(t,L), \dots, \partial_{x_d} F_1(t,L) \right)^\top \mathscr{P}_{\mathcal{N},+} f\] is well defined and
  \[
  \sum_{i=0}^d \| Z^+_i \|^2_{\mathrm{H}^1(\mathrm{U}_\tau \times \mu)} \leq C^{\mathcal N}_{\tau} \|f\|_{\mathrm{L}^2(\mu)}^2,
  \]
  with
  % \|F_0\|^2 + \|\partial_t F_0\|^2 + \sum_i \|\partia_i F_0\|^2 \leq 1/c_\phi + 65 
  % \sum_i \|\partial_i F_1\|^2 + \|\partial_t \partial_i F_1\|^2 + \sum_j \|\partia_i \partial_j F_1\|^2 \leq 36/c_\phi + 72 (...) + 36/c_\phi (...) 
  \begin{equation}
    \label{eq:C_tau_N}
    \begin{aligned}
      C_{\tau}^\mathcal{N} & = \frac{37}{c_\phi} + 65 +  36 \, \left[ 2 + 2 c'_\phi\left(\sqrt{d} + 2 \max \left(8 c'_\phi , \sqrt{c''_\phi \,d} \right)  \right) + \left(1 + \frac{1}{1-\mathrm{e}^{-\tau\,\sqrt{c_\phi}}} \right)^2 \right].  
    \end{aligned}
  \end{equation}

\end{lemma}

\begin{proof}
  Note first that the range of~$\mathscr{P}_{\mathcal{N},+}$ is included by construction in the subspace
  \[
  \mathrm{L}^2_{0,0}(\mathrm{U}_\tau \otimes \mu) = \mathrm{L}^2\left(\mathrm{U}_\tau,\mathrm{L}^2_0(\mu)\right)
  \]
  of functions whose spatial average vanishes for almost all times~$t \in [0,\tau]$ (see the discussion after~\eqref{eq:decomposition_elements_N}). The operators~$F_0(t,L)$ and~$F_1(t,L)$, their time derivatives and their compositions with~$\partial_{x_i}$, make sense on these functional spaces as they involve compositions with inverse powers of~$L$. Since~$\mathscr{P}_{\mathcal{N},+}$ is a projector, the desired bound follows from bounds on various the operators considered on~$\mathcal{B}(\mathrm{L}^2_{0,0}(\mathrm{U}_\tau \otimes \mu))$.
  
  Let us first write some preliminary estimates. The operators~$T_i = \partial_{x_i} L^{-1}$ are bounded on~$\mathrm{L}^2_0(\mu)$ and have an operator norm smaller than~1 since $T_1^\star T_1+\dots+T_d^\star T_d = 1$. Moreover, the operators~$\partial_{x_i} \partial_{x_j} L^{-2}$ for $1 \leq i,j \leq d$ are also bounded on~$\mathrm{L}^2_0(\mu)$. In fact, by~\cite[Lemma~3.6]{bernard2020hypocoercivity} (see also~\cite{dolbeault2015hypocoercivity} and~\cite[Lemma~2.3]{cao2023explicit}),
  \begin{equation}
    \label{eq:eqnew2}
    \sum_{i,j=1}^d \, \left\|\partial_{x_i} \partial_{x_j} L^{-2}\right\|^2_{\mathcal{B}(\mathrm{L}^2_0(\mu))} \leq 2 + 2 c'_\phi\, \left(\sqrt{d} + 2 \max \left(8 c'_\phi , \sqrt{c''_\phi \,d} \right)  \,  \right).
  \end{equation}
  Finally,  
\begin{equation}
  \label{eq:eqnew}
  \sum_{i=1}^d \, \left\|\partial_{x_i}  L^{-2}\right\|^2_{\mathcal{B}(\mathrm{L}^2_0(\mu))} = \left\| L^{-2}\right\|_{\mathcal{B}(\mathrm{L}^2_0(\mu))} \leq \frac{1}{c_\phi}.
\end{equation}

We can then estimate~$F_0(t,L)$ and its derivatives. Rewriting
\[
F_0(t,\zeta) = - \frac{2 \mathrm{e}^{-t\zeta} -1- \mathrm{e}^{-\tau \zeta}}{\zeta} \frac{1-\rme^{-t \zeta}}{1-\rme^{-\tau \zeta}}\frac{1-\rme^{-(\tau-t) \zeta}}{1-\rme^{-\tau \zeta}},
\]
the inequality
% H(\zeta) = 2 \mathrm{e}^{-t\zeta} -1- \mathrm{e}^{-\tau \zeta} is such that H(0) \geq -1, H(1) \leq 1 and H'(\zeta) = 2 e^{-t \zeta} \geq 0
\[
-1 \leq 2 \mathrm{e}^{-t\zeta} -1- \mathrm{e}^{-\tau \zeta} \leq 1,
\]
%the second factor can indeed be bounded as
%\[
%\begin{aligned}
%& \rme^{-t \zeta} - \frac{2}{1+\rme^{-\tau \zeta}}\rme^{-\tau \zeta} \leq \rme^{-t \zeta} - \frac{2}{1+\rme^{\tau \zeta}} \leq \rme^{-t \zeta} - 1 \leq 0, \\
%& \rme^{-t \zeta} - \frac{2}{1+\rme^{-\tau \zeta}}\rme^{-\tau \zeta} \geq \rme^{-\tau \zeta} - \frac{2}{1+\rme^{-\tau \zeta}}\rme^{-\tau \zeta} = \rme^{-\tau \zeta}\frac{\rme^{-\tau \zeta}-1}{1+\rme^{-\tau \zeta}} \geq -1
%\end{aligned}
%\]
immediately implies the bound~$|F_0(t,\zeta)| \leq 1/\zeta$ for any~$t \in [0,\tau]$, so that
\begin{align*}
  \|F_0\|_{\mathcal{B}(\mathrm{L}_{0,0}^2(\mathrm{U}_\tau \otimes \mu))}
  & = \sqrt{ \sup_{\varphi \in \mathrm{L}_{0,0}^2(\mathrm{U}_\tau \otimes \mu)} \frac{1}{\|\varphi\|_{\mathrm{L}^2(\mathrm{U}_\tau \otimes \mu)}^2}\int_0^\tau \|F_0(t,L)\varphi(t,\cdot)\|_{\mathrm{L}^2(\mu)}^2 \, \mathrm{U}_\tau(dt) } \\ 
  & \leq \sup_{t \in [0,\tau]} \|F_0(t,L)\|_{\mathcal{B}(\mathrm{L}_0^2(\mu))} \leq \left\|L^{-1}\right\|_{\mathcal{B}(\mathrm{L}^2_0(\mu))} \leq c_\phi^{-1/2}.
\end{align*}
%Next, 
%\begin{align*}
%  \left| \partial_t F_0(t,\zeta) \right|
%  & \leq \left| \left(1+\rme^{-\tau \zeta}\right) \, \, \rme^{-t \zeta}\,\frac{1-\rme^{-t \zeta}}{1-\rme^{-\tau \zeta}}\frac{1-\rme^{-(\tau-t) \zeta}}{1-\rme^{-\tau \zeta}} \right|  \\
%  & \quad + \left| \left(1+\rme^{-\tau \zeta}\right) \, \left(\rme^{-t \zeta} - \frac{2}{1+\rme^{-\tau \zeta}}\rme^{-\tau \zeta}\right)\frac{\rme^{-t \zeta}}{1-\rme^{-\tau \zeta}}\frac{1-\rme^{-(\tau-t) \zeta}}{1-\rme^{-\tau \zeta}} \right| \\
%  & \quad + \left| \left(1+\rme^{-\tau \zeta}\right) \, \left(\rme^{-t \zeta} - \frac{2}{1+\rme^{-\tau \zeta}}\rme^{-\tau \zeta}\right)\frac{1-\rme^{-t \zeta}}{1-\rme^{-\tau \zeta}}\frac{\rme^{-(\tau-t) \zeta}}{1-\rme^{-\tau \zeta}} \right| \\
%  & \leq 2 \, + \frac{4}{1- \mathrm{e}^{-\tau \, \zeta}},
%\end{align*}
%so that, proceeding as above,
%\[
%\|\partial_t F_0\|_{\mathcal{B}(\mathrm{L}^2_{0,0}(\mathrm{U}_\tau \otimes \mu))} \leq 2 \left(1+\frac{2}{1-\rme^{-\tau \sqrt{c_\phi}}}\right) .
%\]
In order to estimate $\|\partial_{x_i} F_0\|_{\mathcal{B}(\mathrm{L}^2_{0,0}(\mathrm{U}_\tau \otimes \mu))}$ for~$1 \leq i \leq d$, we write~$\partial_{x_i} F_0(t,L) = T_i L F_0(t,L)$ and use
\[
\forall t \in [0,\tau], \qquad \left|\zeta F_0(t,\zeta)\right| \leq \left|\left(2 \mathrm{e}^{-t\zeta} -1- \mathrm{e}^{-\tau \zeta}\right)\frac{1-\rme^{-t \zeta}}{1-\rme^{-\tau \zeta}}\frac{1-\rme^{-(\tau-t) \zeta}}{1-\rme^{-\tau \zeta}}\right| \leq 1,
\]
to conclude that
\[
\sum_{i=1}^d \left\| \partial_{x_i} F_0 \right\|^2_{\mathcal{B}(\mathrm{L}^2_{0,0}(\mathrm{U}_\tau \otimes \mu))} \leq \sum_{i=1}^d \|T_i\|^2_{\mathcal{B}(\mathrm{L}^2_0(\mu))} \sup_{t \in [0,\tau]} \|L F_0(t,L)\|_{\mathcal{B}(\mathrm{L}^2_{0}(\mu))}^2 \leq 1.
\]
Next, using \eqref{eq:div_eq_Z} and~\eqref{eq:from_F_to_mathcal_F}, which imply that~$\partial_t F_0(t,L) =  L F_0(t,L) + L^2 F_1(t,L) - \Id$, it follows that
\begin{align*}
    \|\partial_t F_0\|_{\mathcal{B}(\mathrm{L}^2_{0,0}(\mathrm{U}_\tau \otimes \mu))} \leq 8,
\end{align*}
since~$\|L^2 F_1\|_{\mathcal{B}(\mathrm{L}^2_{0,0}(\mathrm{U}_\tau \otimes \mu))} \leq 6$ by~\eqref{eq:bound_zeta2_F_1} below.

Consider now, for some given~$1 \leq i \leq d$, the operator~$\partial_{x_i} F_1(t,L)$ and its space-time derivatives. Note first that 
\begin{equation}
  \label{eq:bound_zeta2_F_1}
  \forall t \in [0,\tau], \quad \forall \zeta \geq 0, \qquad
  0 \leq \zeta^2 F_1(t,\zeta) = 6 \, \frac{1-\rme^{-t \zeta}}{1-\rme^{-\tau \zeta}}\frac{1-\rme^{-(\tau-t) \zeta}}{1-\rme^{-\tau \zeta}} \mathrm{e}^{-t \zeta} \leq 6. 
\end{equation}
Writing~$\partial_{x_i} F_1(t,L) = \partial_{x_i} L^{-2} L^2 F_1(t,L)$ and combining the above bound with~\eqref{eq:eqnew} gives 
\[
\sum_{i=1}^d \left\| \partial_{x_i} F_1 \right\|^2_{\mathcal{B}(\mathrm{L}^2_{0,0}(\mathrm{U}_\tau \otimes \mu))} \leq \left(\sum_{i=1}^d \left\| \partial_{x_i} L^{-2} \right\|^2_{\mathcal{B}(\mathrm{L}^2_{0}(\mu))}\right)\sup_{t \in [0,\tau]} \left\|L^2F_1(t,L)\right\|^2_{\mathcal{B}(\mathrm{L}^2_{0}(\mu))} \leq \frac{36}{c_\phi}.
\]
A similar decomposition gives, with~\eqref{eq:eqnew2}, 
\begin{align*}
    \sum_{i,j=1}^d \left\| \partial_{x_i} \partial_{x_j} F_1 \right\|^2_{\mathcal{B}(\mathrm{L}^2_{0,0}(\mathrm{U}_\tau \otimes \mu))} &\leq 72 \left[ 1 + c'_\phi\left(\sqrt{d} + 2 \max \left(8 c'_\phi , \sqrt{c''_\phi \,d} \right)  \right) \right].
\end{align*}
Finally, for~$\zeta \geq \sqrt{c_\phi}$, the bound
\[
\begin{aligned}
  \left| \zeta \partial_t F_1(t,\zeta) \right| & \leq 6 + 6 \left| \frac{\mathrm{e}^{-t\zeta}}{1-\mathrm{e}^{-\tau \zeta}} \, \frac{1-\mathrm{e}^{-(\tau-t)\zeta}}{1-\mathrm{e}^{-\tau \zeta}} - \frac{1-\mathrm{e}^{-t \zeta}}{1-\mathrm{e}^{-\tau \zeta}} \, \frac{\mathrm{e}^{-(\tau-t)\zeta}}{1-\mathrm{e}^{-\tau \zeta}} \right| \\
  & = 6 + 6 \frac{\left| \mathrm{e}^{-t\zeta} -\mathrm{e}^{-(\tau-t) \zeta}\right|}{(1-\mathrm{e}^{-\tau \zeta})^2} 
  \leq 6 \left( 1 + \frac{1}{1-\mathrm{e}^{-\tau \sqrt{c_\phi}}}\right),
\end{aligned}
\]
and the decomposition~$\partial_t \partial_{x_i} F_1(t,L) = \partial_{x_i} L^{-1} L \partial_t F_1(t,L)$ imply that
\[
\begin{aligned}
  \sum_{i=1}^d \left\| \partial_t \partial_{x_i} F_1 \right\|_{\mathcal{B}(\mathrm{L}^2_{0,0}(\mathrm{U}_\tau \otimes \mu))}^2 & \leq \left( \sum_{i=1}^d \left\| \partial_{x_i} L^{-1} \right\|^2_{\mathcal{B}(\mathrm{L}^2_0(\mu))} \right) \sup_{t \in [0,\tau]} \left\|L \partial_t F_1(t,L)\right\|^2_{\mathcal{B}(\mathrm{L}^2_{0}(\mu))} \\
  & \leq 36 \left(1 + \frac{1}{1-\mathrm{e}^{-\tau\,\sqrt{c_\phi}}} \right)^2. 
\end{aligned}
\]
Gathering the estimates allows to conclude the proof.
\end{proof}

We are now in position to conclude the proof of Proposition~\ref{prop:estimates_divergence_equation}.

\begin{proof}[Proof of Proposition~\ref{prop:estimates_divergence_equation}]
  We decompose the source term $f$ of \eqref{eq:Lions_A} into its component along the orthogonal subspaces~$\mathcal{N}^\perp$ and~$\mathcal{N}$. Since~\eqref{eq:Lions_A} is linear, the solution is obtained by summing the solutions with right hand side $\mathscr{P}_{\mathcal{N}} f$ and $\mathscr{P}_{\mathcal{N}^\perp}f$, as given by~\eqref{eq:explicit_solution_divergence_equation}. For the part of the solution associated with the component in~$\mathcal{N}^\perp,$ we apply Lemma~\ref{lem:Lax_Milgram}; while for the part of the solution associated with the component in~$\mathcal{N}$ we use Lemma~\ref{Nlem}. The constant~$C_\tau^{\mathrm{div}}$ in~\eqref{eq:estimates_solution_divergence_equation} is finally bounded as
  \begin{equation}
    \label{eq:final_constant_Lions_ineq}
    \left(C_\tau^\mathrm{div}\right)^2 \leq 3\left(C_\tau^\mathrm{LM}+2C_\tau^\mathcal{N}\right),
  \end{equation}
  these constants being respectively defined in~\eqref{eq:C_tau_LM} and~\eqref{eq:C_tau_N}.
  The terms in~\eqref{eq:C_tau_LM} can be bounded by~$2 + C_\phi \sqrt{d}\tau^2 $, while the terms in~\eqref{eq:C_tau_N} can be bounded by~$C_\phi(\sqrt{d}+\tau^{-2})$, for an explicit constant $C_\phi$ which does not directly depend on the dimension $d$, but only on~$c_\phi,c'_\phi,$ and~$c''_\phi$.
\end{proof}

\begin{remark}
  \label{rmk:scaling_Lions_constant_Poincare}
  To obtain the scaling of~$C_\tau^\mathrm{Lions}$ with respect to~$c_\phi$ in the regime~$c_\phi \to 0$, we rely on the upper bound given by~\eqref{eq:final_constant_Lions_ineq}. In view of~\eqref{eq:C_tau_LM}, the constant~$C_\tau^\mathrm{LM}$ scales as~$1/c_\phi$ when choosing~$\tau$ smaller than~$c_\phi^{-1/2}$. Moreover, \eqref{eq:C_tau_N} shows that the constant~$C_\tau^\mathcal{N}$ scales as~$1/c_\phi$ for~$\tau$ larger than~$c_\phi^{-1/2}$. Overall, the choice~$\tau \sim c_\phi^{-1/2}$ therefore implies that~$C_\tau^\mathrm{Lions}$ scales as~$1/c_\phi$. More precisely,
  \begin{equation}
      \label{eq:R_Lions}
      C_\tau^\mathrm{Lions} = \frac{1}{c_\phi} \left[ C_1\left(c_\phi',c_\phi''\right) \left(c_\phi+\sqrt{d}\right) + C_2 \right],
  \end{equation}
  for some constants~$C_1,C_2$, with~$C_1$ depending on~$\phi$ only through~$c_\phi',c_\phi''$. 
\end{remark}

%-------------------------
\section{Scaling of the exponential convergence rate for tensorised kinetic energies}
\label{sec:tenso}

We state and prove in this section a result on the dimension dependence of the constant~$K_\mathrm{avg}$ of Lemma~\ref{avglem}, for tensorised kinetic energies (see the discussion in Section~\ref{sec:scaling_dimension}).

\begin{proposition}\label{kavgtens}
  Suppose that the kinetic energy is of the separable form~\eqref{tensorised} and satisfies the conditions~\eqref{eq:integrability_conditions_on_psi2}, \eqref{eq:coercivity_Delta_psi_OU2} and~\eqref{eq:condition_vanishing}. Then, the constant~$K_{\mathrm{avg}}$ of Lemma~\ref{avglem} can be chosen uniformly with respect to the dimension~$d$. More precisely,
  \begin{equation}
  \label{eq:K_avg_q}
  K_{\mathrm{avg}}^q = 2\max \left(\mathscr{C}_1,\mathscr{C}_2\right)^2,
  \end{equation}
  with
  \[
\begin{aligned}
  \mathscr{C}_1 & = 1 + \frac{\|q'\|_{\mathrm{H}^1(\mathrm{e}^{-q})}}{\|q'\|_{\mathrm{L}^2(\mathrm{e}^{-q})}^2},\\
  \mathscr{C}_2 & = 1 + \|q'\|_{\mathrm{L}^2(\mathrm{e}^{-q})} + \frac{1}{\|q'\|_{\mathrm{L}^2(\mathrm{e}^{-q})}} + \frac{\left\|(q')^2-q''\right\|_{\mathrm{L}^2(\mathrm{e}^{-q})} + \left\|q''\right\|_{\mathrm{L}^2(\mathrm{e}^{-q})}}{\left\|q'\right\|_{\mathrm{L}^2(\mathrm{e}^{-q})}^2} .
\end{aligned}
\]
\end{proposition}

\begin{proof}
  We refine the proof of Lemma \ref{avglem}. Note that the matrices defined in~\eqref{eq:matrix_M} and~\eqref{eq:matrix_M_2} are respectively 
  \[
  \mathscr{M} = \|q'\|^2_{\mathrm{L}^2(\mathrm{e}^{-q})} \mathrm{Id}_d, \qquad M = d \mathrm{Id}_d,
  \]
  so that~$\rho(\mathscr{M}) = \|q'\|_{\mathrm{L}^2(\mathrm{e}^{-q})}^2$ and~$\rho(M) = d$. The estimate~\eqref{eq:dt_Pi_h_H-1} on the time derivative can therefore be used as such. On the other hand, the estimates on the spatial derivatives need to be made more precise. Recall that we consider a test function~$Z=(Z_1,\dots,Z_d) \in \mathrm{H}^1_\ZDCT(\mathrm{U}_\tau \otimes \mu)^d$ with
  \[
  \|Z\|^2_{\mathrm{H}^1(\mathrm{U}_\tau \otimes \mu)} = \sum_{i=1}^d \|Z_i\|^2_{\mathrm{L}^2(\mathrm{U}_\tau \otimes \mu)} + \|\partial_t Z_i\|^2_{\mathrm{L}^2(\mathrm{U}_\tau \otimes \mu)} + \sum_{i,j = 1}^d \left\|\partial_{x_j} Z_i \right\|^2_{\mathrm{L}^2(\mathrm{U}_\tau \otimes \mu)} \leq 1. 
  \]
  We restart from~\eqref{eq:series_eq_nabla_x_Pi_Z}, with~$G_i(v) = d^{-1/2} q'(v_i)/\|q'\|_{\mathrm{L}^2(\mathrm{e}^{-q})}$. The factor~$\|\nabla_v \psi\|_{\mathrm{L}^2(\gamma)}^{-1} \|G \cdot MZ\|_{\mathrm{L}^2(\mathrm{U}_\tau \otimes \mu,\mathrm{H}^1(\gamma))}$ arising when estimating the first integral in the last equality of~\eqref{eq:series_eq_nabla_x_Pi_Z} can be bounded by writing
  \[
  \begin{aligned}
    & \frac{1}{\|\nabla_v \psi\|_{\mathrm{L}^2(\gamma)}^2} \|G \cdot MZ\|_{\mathrm{L}^2(\mathrm{U}_\tau \otimes \mu,\mathrm{H}^1(\gamma))}^2
    = \frac{1}{\|q'\|_{\mathrm{L}^2(\mathrm{e}^{-q})}^4} \left\| \sum_{i=1}^d q'(v_i) Z_i \right\|_{\mathrm{L}^2(\mathrm{U}_\tau \otimes \mu,\mathrm{H}^1(\gamma))}^2 \\
    & \qquad \leq \frac{1}{\|q'\|_{\mathrm{L}^2(\mathrm{e}^{-q})}^4}\left( \left\| \sum_{i=1}^d q'(v_i) Z_i \right\|_{\mathrm{L}^2(\mathrm{U}_\tau \otimes \Theta)}^2 + \sum_{k=1}^d \left\| \partial_{v_k} \left( \sum_{i=1}^d q'(v_i) Z_i \right) \right\|_{\mathrm{L}^2(\mathrm{U}_\tau \otimes \Theta)}^2 \right). 
  \end{aligned}
  \]
  Now, using that~$v_i \mapsto q'(v_i)$ has average~0 with respect to~$\gamma$,
  \[
  \begin{aligned}
    \left\| \sum_{i=1}^d q'(v_i) Z_i \right\|_{\mathrm{L}^2(\mathrm{U}_\tau \otimes \Theta)}^2 & = \sum_{i,j=1}^d \left(\int_\mathcal{V} q'(v_i) q'(v_j) \, \gamma(dv)\right)\left(\int_0^\tau \int_\mathcal{X} Z_i Z_j \, d\mathrm{U}_\tau \, d\mu\right) \\
    & = \left\|q'\right\|_{\mathrm{L}^2(\mathrm{e}^{-q})}^2 \sum_{i=1}^d \|Z_i\|^2_{\mathrm{L}^2(\mathrm{U}_\tau \otimes \mu)},
  \end{aligned}
  \]
  and
  \[
  \left\| \partial_{v_k} \left( \sum_{i=1}^d q'(v_i) Z_i \right) \right\|_{\mathrm{L}^2(\mathrm{U}_\tau \otimes \Theta)}^2 = \left\| q''(v_k) Z_k \right\|_{\mathrm{L}^2(\mathrm{U}_\tau \otimes \Theta)}^2 = \left\|q''\right\|_{\mathrm{L}^2(\mathrm{e}^{-q})}^2 \|Z_k\|^2_{\mathrm{L}^2(\mathrm{U}_\tau \otimes \mu)},
  \]
  so that
  \[
  \frac{1}{\|\nabla_v \psi\|_{\mathrm{L}^2(\gamma)}^2} \|G \cdot MZ\|_{\mathrm{L}^2(\mathrm{U}_\tau \otimes \mu,\mathrm{H}^1(\gamma))}^2 \leq \frac{\|q'\|_{\mathrm{H}^1(\mathrm{e}^{-q})}^2}{\|q'\|_{\mathrm{L}^2(\mathrm{e}^{-q})}^4}.
  \]
  For the third integral in the last equality of~\eqref{eq:series_eq_nabla_x_Pi_Z}, the factor involving~$Z$ is bounded as
  \[
  \frac{1}{\|\nabla_v \psi\|_{\mathrm{L}^2(\gamma)}} \|MG \cdot \partial_t Z\|_{\mathrm{L}^2(\mathrm{U}_\tau \otimes \Theta)} \leq \frac{1}{\|q'\|_{\mathrm{L}^2(\mathrm{e}^{-q})}},  
  \]
  by estimates similar to the previous ones, upon replacing~$Z$ by~$\partial_t Z$ and considering a~$\mathrm{L}^2(\gamma)$ instead of a $\mathrm{H}^1(\gamma)$ norm. 

  It remains at this stage to estimate the second integral in the last equality of~\eqref{eq:series_eq_nabla_x_Pi_Z}, which follows by estimating
  \begin{equation}
    \label{eq:two_terms_T_to_estimate}
    \frac{\|T(MG \cdot Z)\|_{\mathrm{L}^2(\mathrm{U}_\tau \otimes \Theta)}}{\|\nabla_v \psi\|_{\mathrm{L}^2(\gamma)}} = \frac{1}{\|q'\|_{\mathrm{L}^2(\mathrm{e}^{-q})}^2}\left\| T\left(\sum_{i=1}^d q'(v_i) Z_i \right)\right\|_{\mathrm{L}^2(\mathrm{U}_\tau \otimes \Theta)}.  
  \end{equation}
  Since
  \[
  T = \sum_{i=1}^d \partial_{v_i}^\star \partial_{x_i} - \partial_{v_i} \partial_{x_i}^\star,
  \]
  and denoting by~$Q(v_i) = \partial_{v_i}^\star q'(v_i) = \left(\partial_{v_i}^\star\right)^2 \mathbf{1} = q'(v_i)^2 - q''(v_i)$, it follows that
  \[
  T\left(\sum_{i=1}^d q'(v_i) Z_i \right) = \sum_{i=1}^d Q(v_i) \partial_{x_i} Z_i - q''(v_i) \partial_{x_i}^\star Z_i.
  \]
  Since~$Q$ has average~0 with respect to~$\gamma$, and since~\eqref{eq:condition_vanishing} implies
  \[
  \int_\mathbb{R} Q q'' \, \rme^{-q} = \int_\mathbb{R} \left[ \left(\partial_{v_1}^\ast\right)^2 \mathbf{1}\right] \partial_{v_1}^2 q(v_1) \, \mathrm{e}^{-q(v_1)} dv_1 = \int_\mathbb{R} q^{(4)} \, \rme^{-q} = 0,
  \]
  we obtain that
  \begin{align}
    & \left\| T\left(\sum_{i=1}^d q'(v_i) Z_i \right)\right\|_{\mathrm{L}^2(\mathrm{U}_\tau \otimes \Theta)}^2 = \int_0^\tau \int_\mathcal{X} \int_\mathcal{V}  \left( \sum_{i=1}^d Q \partial_{x_i} Z_i - q'' \partial_{x_i}^\star Z_i \right)^2 d\mathrm{U}_\tau \, d\Theta \notag\\
    & = \|Q\|_{\mathrm{L}^2(\mathrm{e}^{-q})}^2 \sum_{i=1}^d\left\| \partial_{x_i} Z_i \right\|_{\mathrm{L}^2(\mathrm{U}_\tau \otimes \mu)}^2 + \sum_{i,j = 1}^d \int_\mathcal{V} q''(v_i) q''(v_j) \, \gamma(dv) \int_0^\tau \int_\mathcal{X} \partial_{x_i}^\star Z_i \partial_{x_j}^\star Z_j \, d\mathrm{U}_\tau \, d\mu \label{eq:T_terms_to_control_for_Kavg} \\
    & = \|Q\|_{\mathrm{L}^2(\mathrm{e}^{-q})}^2 \sum_{i=1}^d\left\| \partial_{x_i} Z_i \right\|_{\mathrm{L}^2(\mathrm{U}_\tau \otimes \mu)}^2 + \sum_{i,j = 1}^d \int_\mathcal{V} q''(v_i) q''(v_j) \, \gamma(dv) \int_0^\tau \int_\mathcal{X} \partial_{x_j} Z_i \partial_{x_i} Z_j \, d\mathrm{U}_\tau \, d\mu \notag \\
    & \leq \|Q\|_{\mathrm{L}^2(\mathrm{e}^{-q})}^2 \sum_{i=1}^d\left\| \partial_{x_i} Z_i \right\|_{\mathrm{L}^2(\mathrm{U}_\tau \otimes \mu)}^2 + \sum_{i,j = 1}^d \int_\mathcal{V} \left|q''(v_i) q''(v_j) \right| \gamma(dv) \frac{\left\| \partial_{x_i} Z_j \right\|_{\mathrm{L}^2(\mathrm{U}_\tau \otimes \mu)}^2 + \left\| \partial_{x_j} Z_i \right\|_{\mathrm{L}^2(\mathrm{U}_\tau \otimes \mu)}^2}{2}. \notag
  \end{align}
  The integral in~$v$ in the second term of the last inequality is bounded by~$\left\|q''\right\|_{\mathrm{L}^2(\mathrm{e}^{-q})}^2$ (using a Cauchy--Schwarz inequality for~$i \neq j$), so that 
  \[
  \begin{aligned}
    \left\| T\left(\sum_{i=1}^d q'(v_i) Z_i \right)\right\|_{\mathrm{L}^2(\mathrm{U}_\tau \otimes \Theta)}^2 & \leq \|Q\|_{\mathrm{L}^2(\mathrm{e}^{-q})}^2 \sum_{i=1}^d \left\| \partial_{x_i} Z_i \right\|_{\mathrm{L}^2(\mathrm{U}_\tau \otimes \mu)}^2 + \left\|q''\right\|_{\mathrm{L}^2(\mathrm{e}^{-q})}^2 \sum_{i,j=1}^d \left\| \partial_{x_i} Z_j \right\|_{\mathrm{L}^2(\mathrm{U}_\tau \otimes \mu)}^2 \\
    & \leq \left( \|Q\|_{\mathrm{L}^2(\mathrm{e}^{-q})}^2 + \left\|q''\right\|_{\mathrm{L}^2(\mathrm{e}^{-q})}^2 \right) \left\| Z \right\|_{\mathrm{H}^1(\mathrm{U}_\tau \otimes \mu)}^2, 
  \end{aligned}
  \]
  and finally
  \[
  \left\| T\left(\sum_{i=1}^d q'(v_i) Z_i \right)\right\|_{\mathrm{L}^2(\mathrm{U}_\tau \otimes \Theta)} \leq \|Q\|_{\mathrm{L}^2(\mathrm{e}^{-q})} + \left\|q''\right\|_{\mathrm{L}^2(\mathrm{e}^{-q})}.
  \]
  We obtain the sought expression for~$K_{\mathrm{avg}}$ by gathering the above estimates in view of~\eqref{kavg}.
\end{proof}

\begin{remark}
\label{rmk:dim_dep_Kavg}
  When~\eqref{eq:condition_vanishing} is not satisfied, and since~$Q$ has average~0 with respect to~$\gamma$, there is an extra term
  \[
  -2 \left( \int_\mathcal{V} Q q'' \, \rme^{-q}\right) \sum_{i=1}^d \int_0^\tau \int_\mathcal{X} \partial_{x_i} Z_i \partial_{x_i}^\star Z_i \, d\mathrm{U}_\tau \, d\mu
  \]
  to be added to~\eqref{eq:T_terms_to_control_for_Kavg}. The space-time integral can be rewritten as 
  \[
  \int_0^\tau \int_\mathcal{X} \partial_{x_i} Z_i \partial_{x_i}^\star Z_i \, d\mathrm{U}_\tau \, d\mu = - \left\| \partial_{x_i} Z_i \right\|_{\mathrm{L}^2(\mathrm{U}_\tau \otimes \Theta)}^2 + \int_0^\tau \int_\mathcal{X} \left(\partial_{x_i} Z_i\right) \left(\partial_{x_i} \phi\right) Z_i \, d\mathrm{U}_\tau \, d\mu.
  \]
  The first term on the right hand side can be controlled in absolute value by~$\left\| Z \right\|_{\mathrm{H}^1(\mathrm{U}_\tau \otimes \mu)}^2 \leq 1$ when summing over~$1 \leq i \leq d$. On the other hand, the second term cannot be controlled uniformly in the dimension~$d$. In fact, it typically scales as~$\sqrt{d}$, as the following computations show when the potential energy function is of the separable form~\eqref{eq:separable_pot}. First, an upper bounded of order~$\sqrt{d}$ is obtained via a Cauchy--Schwarz inequality: 
  \[
  \left| \sum_{i=1}^d \int_0^\tau \int_\mathcal{X} \partial_{x_i} Z_i \partial_{x_i} \phi \, d\mathrm{U}_\tau \, d\mu \right| \leq \left\| \nabla \phi \right\|_{\mathrm{L}^2(\mu)}\sqrt{\sum_{i=1}^d \left\| \partial_{x_i} Z_i \right\|_{\mathrm{L}^2(\mathrm{U}_\tau \otimes \mu)}^2} \leq \left\| \nabla \phi \right\|_{\mathrm{L}^2(\mu)},
  \]
  the latter quantity being equal to~$\sqrt{d} \| \overline{\phi}' \|_{\mathrm{L}^2(\mathrm{e}^{-\bar{\phi}})}$ for separable potentials. Moreover, choosing~$Z_i(t,x) = d^{-1/2} \bar{\phi}(x_i)/\left\| \bar{\phi} \right\|_{\mathrm{H}^1(\mathrm{e}^{-\bar{\phi}})}$ for~$1 \leq i \leq d$ ensures that~$\|Z\|_{\mathrm{H}^1(\mathrm{U}_\tau \otimes \mu)} = 1$ and
  \[
  \sum_{i=1}^d  \int_0^\tau \int_\mathcal{X} \partial_{x_i} Z_i \partial_{x_i} \phi \, d\mathrm{U}_\tau \, d\mu= \sqrt{d} \frac{\left\| \overline{\phi}' \right\|_{\mathrm{L}^2(\mathrm{e}^{-\bar{\phi}})}^2}{\left\| \overline{\phi}' \right\|_{\mathrm{H}^1(\mathrm{e}^{-\bar{\phi}})}}.
  \]
  Therefore, the constant~$K_\mathrm{avg}$ scales as~$\sqrt{d}$ when~\eqref{eq:condition_vanishing} is not satisfied, instead of being bounded uniformly with respect to the dimension as in Proposition~\ref{kavgtens}.
\end{remark}

\section*{Acknowledgements}
The authors thank Emeric Bouin, Evan Camrud, Jean Dolbeault, Susanne Pflügl, and Lihan Wang for stimulating discussions and feedbacks. 

G.B.~has been funded by the European Union’s Horizon 2020 research and innovation program under the Marie Sklodow\-ska-Curie grant agreement No.~754362, and by the Project EFI (ANR-17-CE40-0030) of the French National Research Agency (ANR). The work of G.S. was partially funded by the European Research Council (ERC) under the European Union's Horizon 2020 research and innovation programme (grant agreement No~810367), and by the Agence Nationale de la Recherche under grants ANR-19-CE40-0010 (QuAMProcs) and ANR-21-CE40-0006 (SINEQ).

%--------------------------------------
\appendix
\section{Sufficient conditions for hypoellipticity}
\label{app:Hormander}

Regularity results can be justified thanks to the hypoelliptic regularity of~\cite{hormander1967hypoelliptic}. Indeed, the Lie algebra generated by $\nabla_v$ and $T$ is of dimension $2d$ at all points~$(x,v) \in \mathcal{X} \times \mathcal{V}$ when Assumption~\ref{ass3} holds. A sufficient condition to this end is that the Hessian~$\nabla^2_{vv} \psi$ is of full rank for all~$v \in \mathcal{V}$. Indeed, the commutator between the Hamiltonian part~$T$ and the velocity gradient is
\[
  [T,\nabla_v] = \nabla_v \psi \cdot \nabla_x \, \nabla_v - \nabla_x \phi \cdot \nabla_v \, \nabla_v - \nabla_v ( \nabla_v \psi \cdot \nabla_x) + \nabla_v \, \nabla_x \phi \cdot \nabla_v = \left( -\nabla^2_{vv} \psi \right) \nabla_x,
  \]
  so that all derivatives in~$x$ can be recovered using this commutator. If $\nabla^2_{vv} \psi$ is not of full rank somewhere, further commutators should be computed, in the spirit of \cite{villani2009hypocoercivity}, which will involve higher order derivatives of~$\psi$. For instance,
  \[
  \left[\partial_{v_j}, \left[T,\partial_{v_i}\right] \right] = -\left[\partial_{v_j}, \nabla_v \left(\partial_{v_i} \psi\right) \cdot \nabla_x \right] = -\nabla_v \left(\partial^2_{v_i,v_j} \psi\right) \cdot \nabla_x,
  \]
  so that hypoellipticity is ensured when~$\{ \nabla_v (\partial^2_{v_i,v_j} \psi) \}_{1 \leq i,j \leq d}$ has full rank~$d$ at those~$v \in \mathcal{V}$ for which~$\nabla^2_{vv} \psi(v)$ is not of full rank.

%----------------

\bibliographystyle{siam}
\bibliography{biblio,Newref}
\end{document}